\documentclass[reqno]{amsart}
 \usepackage{amsbsy,amssymb,amscd,amsfonts,latexsym,amstext,delarray,
 amsmath,graphicx, color}
 \usepackage[dvips]{epsfig}
\usepackage{hyperref}
\input xypic

\definecolor{Green}{rgb}{0,1,0}
\definecolor{Blue}{RGB}{0,0,191}
\definecolor{mathmodecolor}{RGB}{0,102,0}
\definecolor{keywordcolor}{RGB}{0,51,151}
\definecolor{sourcebackgroundcolor}{RGB}{255,247,223}
\definecolor{unixagred}{RGB}{255,0,0}
\definecolor{lightgray}{RGB}{191,191,191}
\definecolor{green}{RGB}{1,191,191}

\newtheorem{thm}{Theorem}[section]
\newtheorem{prop}[thm]{Proposition}
\newtheorem{cor}[thm]{Corollary}
\newtheorem{lem}[thm]{Lemma}

\newtheorem{defn}[thm]{Definition}
\newtheorem{rem}[thm]{Remark}
\newtheorem{example}[thm]{Example}

\def\qqq{\,,\quad~\forall}

\def\Hom{{\rm Hom}}

\def\Ker{{\rm Ker}}

\def\\adsod{{\rm \adsod}}

\def\sign{{\rm sign}}

\def\Spec{{\rm Spec\,}}
\def\Sp{{\rm Spec}\,}

\def\A{{\mathbb A}}
\def\B{{\mathbb B}}
\def\C{{\mathbb C}}
\def\F{{\mathbb F}}
\def\K{{\mathbb K}}

\def\N{{\mathbb N}}
\def\P{{\mathbb P}}
\def\Q{{\mathbb Q}}
\def\R{{\mathbb R}}
\def\Z{{\mathbb Z}}
\def\H{{\mathbb H}}

\def\cD{{\mathcal D}}

\def\cH{{\mathcal H}}

\def\cL{{\mathcal L}}
\def\c\ads{{\mathcal \ads}}

\def\cO{{\mathcal O}}
\def\cP{{\mathcal P}}

\def\cW{{\mathcal W}}

\newcommand{\ie}{{\it i.e.\/}\ }
\newcommand{\eg}{{\it e.g.\/}\ }
\newcommand{\cf}{{\it cf.\/}\ }
\newcommand{\opcit}{{\it op.cit.\/}\ }

\def\dim{{\mbox{dim}}}

\def\Hom {{\mbox{Hom}}}

\def\ffp{\mathfrak{p}}

\def\\adso{\mathfrak{\adso}}
\def\An{\mathfrak{ Ring}}
\def\han{\mathfrak{ Hring}}

\def\\adsr{\mathfrak{ \adsR}}

\def\Se{\mathfrak{ Sets}}

\def\rmax{\R_+^{\rm max}}
\def\ads{\H_\K}
\def\kras{\mathbf{K}}
\def\sign{\mathbf{S}}

\newcommand{\nil}[1]{}

\title
{The hyperring of ad\`ele classes}

\author[Connes]{Alain Connes}
\author[Consani]{Caterina Consani}
\address{A.~Connes: Coll\`ege de France \\
3, rue d'Ulm \\ Paris, F-75005 France
\\ I.H.E.S. and Vanderbilt
University} \email{alain\@@connes.org}
\address{C.~Consani: Mathematics Department \\ Johns Hopkins
University \\ Baltimore, MD 21218 USA} \email{kc\@@math.jhu.edu}
\thanks{The authors are partially supported by the NSF grant
DMS-FRG-0652164. The first author thanks S. Barr\'e, E. Ghys and B. Sevennec for references on non-Desarguesian planes.}


\begin{document}

\maketitle

\begin{abstract} We show that the theory of hyperrings, due to M. Krasner, supplies a perfect framework to understand the algebraic structure of the ad\`ele class space $\H_\K=\A_\K/\K^\times$ of a global field $\K$. After promoting $\F_1$ to a hyperfield $\kras$, we prove that a hyperring of the form $R/G$ (where $R$ is a ring and $G\subset R^\times$  is a subgroup of its multiplicative group) is a hyperring extension of $\kras$ if and only if $G\cup\{0\}$ is a subfield of $R$. This result applies  to the ad\`ele class space  which thus  inherits the structure of a hyperring extension $\H_\K$ of $\kras$. We begin to investigate the content of an algebraic geometry over $\kras$. The category of commutative hyperring extensions of $\kras$ is inclusive of:  commutative algebras over fields with semi-linear homomorphisms,  abelian groups with injective homomorphisms and  a rather exotic land  comprising homogeneous non-Desarguesian planes. Finally, we show that for a global field $\K$ of positive characteristic, the groupoid of the prime elements of the hyperring $\H_\K$ is canonically and equivariantly isomorphic to the groupoid of the loops of the maximal abelian cover of the curve  associated to the global field $\K$.
\end{abstract}

\tableofcontents

\section{Introduction}

The goal of this paper is to understand the algebraic structure of the ad\`ele class space $\H_\K=\A_\K/\K^\times$ of a global field $\K$.
In our recent work \cite{announc3}, we have shown that   the introduction of an elementary theory of algebraic geometry over the absolute point $\Spec\F_1$ reveals the role of the natural monoidal structure of the ad\`ele class space $\A_\K/\K^\times$ of a global field. This structure is used  to reformulate, in a more conceptual manner,  the spectral realization of zeros of $L$-functions. In the subsequent paper \cite{jamifine}, we have given substantial evidence to the statement that idempotent analysis and tropical geometry determine, through the theory of idempotent semi-rings, a natural framework where to develop mathematics in ``characteristic one".

 A key  role in the formulation of these ideas is played by the procedure of de-quantization that requires the replacement of the use of real analysis by its idempotent version, and the implementation of the semifield $\rmax$ in place of the classical $\R_+$. Long ago,  M. Krasner devised an analogous procedure that can be performed at a finite place of $\Q$ (\cf \cite{Krasner}). His construction shows how to approximate a local field $\F_q((T))$ of positive characteristic by a system of local fields of characteristic zero and with the same residue field, as the absolute ramification index tends to infinity. Krasner's method is based on the idea of class field and on the generalization of the classical additive law in a ring  by the structure of a hypergroup in the sense of F. Marty \cite{Marty}. This process produces the notion of a (Krasner) {\em hyperring} (\cf\cite{Krasner1}) which fits perfectly with our previous constructions and in particular with the framework of noncommutative geometry.

In the usual theory of semi-rings, it is not possible to reconcile the characteristic one property stating that $x+x=x$ for all elements $x$ of a semi-ring $R$, with the additive group law requiring that every element in $R$ admits an additive inverse. On the other hand, the existence of an additive inverse plays a crucial role when, for instance, tensor products are involved. The structure of a hyperring makes this compatibility --between characteristic one and existence of additive inverse-- possible. Remarkably, the ad\`ele class space $\ads=\A_\K/\K^\times$ of a global field $\K$  turns out to possess the correct hyperring structure that combines the two above properties and in particular one has
$
x+x=\{0,x\}
$
for all $x\in \ads$.

This formula means that $\ads$ is a hyperring over the simplest hyperfield $\kras$ that is defined as the set $\{0,1\}$ endowed with the obvious multiplication and a hyper-addition requiring that $1+1=\{0,1\}$. Moreover, while the quotient of a ring $R$ by a subgroup $G\subset R^\times$ of its multiplicative group is always a hyperring (\cf\cite{Krasner1}), we find that $R/G$ is an extension of $\kras$ exactly when $G\cup \{0\}$ is a {\em subfield} of $R$ (\cf Proposition \ref{krasner2}).

We explicitly remark here that the ``absolute point" $\Sp \F_1$ should not be confused with $\Sp \kras$, in fact while $\Sp \F_1$  sits under $\Spec\Z$, $\Sp \kras$   is the natural lift of $\Sp\F_1$ above the generic point of $\Spec\Z$.
 \begin{gather}
\label{overall0}
 \,\hspace{100pt}\raisetag{-47pt} \xymatrix@C=25pt@R=25pt{
 \Spec\Z\ar[d] &
  \Spec\kras\ar[dl]\ar[l]& &\\
\Spec\F_1   & &\\
}\hspace{25pt}
\end{gather}
In this paper we show that after suitably extending the classical definition of a $\Z$-scheme, by replacing the category of (commutative) rings with that of hyperrings (as was done \eg in  \cite{Procesi}), the spectrum $\Sp \kras$ plays the role of the ``generic point" in algebraic geometry. In fact, in Proposition \ref{ex2} we prove that for any scheme $X$ of finite type over $\Z$, there is a canonical identification of sets
\begin{equation}\label{identsets}
    X\simeq \Hom(\Sp(\kras),X).
\end{equation}

 One should not confuse the content of a geometry over $\Sp\F_1$, that  essentially means a theory of (pointed) monoids (\cf\cite{deit} and \cite{announc3}), with the more refined geometric theory over $\Sp\kras$ that no longer ignores the {\em additive} structure.
 For instance, one finds that the prime spectrum of the {\em monoid} $\A_\K/\K^\times$ involves all subsets of the set $\Sigma_\K$ of places of the global field $\K$, while  the prime spectrum of the {\em hyperring} $\ads$ is made by the subsets of $\Sigma_\K$ with only {\em one element}. By  restricting this study to the ideals which are {\em closed} in the natural topology, one  obtains the natural identification $\Sp(\ads)=\Sigma_\K$.

The examples of tensor products of hyperrings that we consider in this paper,  which correspond geometrically to the fiber product of the geometric spectra, allow us to understand, at a more conceptual level, several fundamental constructions of noncommutative geometry. In particular, this provides a  new perspective on the structure of the BC-system \cite{ccm}.\vspace{.05in}

The rule of signs is a basic principle in elementary arithmetic. It is a simple fact that while the sign of the product of two numbers is uniquely determined by their respective signs,  the  sign of the sum of a positive and a negative number is ambiguous (\ie it can be $+,-,0$). As a straightforward encoding of this rule, one can upgrade the  monoid $\F_{1^2}$ into a hyperfield with three elements: $\sign=\{-1,0,1\}$. Following this viewpoint, one discovers that the BC-system is directly related to the following hyperring extension of $\sign$
$$
\Z_\sign:=\hat\Z\otimes_\Z\sign,
$$
which is obtained by implementing the natural sign homomorphism $\Z\to \sign$ and the embedding $\Z\to \hat\Z$ of the relative integers into the profinite completion. By taking the topological structure into account, the spectrum $\Sp(\Z_\sign)$ is isomorphic to $\Sp(\Z)$, but unlike this latter space, $\Sp(\Z_\sign)$ maps naturally to $\Sp\sign$. Incidentally,  we remark that the map $\Sp(\Z_\sign)\to \Sp\sign$ should be viewed as a refinement (and a lift) of the obvious map $\Sp(\Z)\to\Sp\F_1$.

The process of adjoining the archimedean place is obtained by moving from finite ad\`eles to the full ad\`eles $\A_\Q$ over $\Q$. Following the hyperring structures, one sees that the hyperfield $\kras$ is the quotient of $\sign$ by the subgroup $\{\pm 1\}$. This fact determines a canonical surjection (absolute value) $\pi: \sign \to \kras$ which is used to show that the ad\`ele class space is described by the hyperring
$$
\H_\Q=\A_\Q\otimes_\Z\kras
$$
whose associated spectrum is  $\Sp(\H_\Q)=\Sp(\Z)\cup\{ \infty\}=\Sigma_\Q$.\vspace{.05in}

In \S \ref{sectproj}, we take the viewpoint of W. Prenowitz   \cite{Prenowitz} and R. Lyndon \cite{Lyndon} to explain a natural correspondence between $\kras$-vector spaces and projective geometries in which every line has at least $4$ points. By implementing some classical results of incidence geometry mainly due to H. Karzel \cite{Karzel}, we describe the classification of finite hyperfield extensions of $\kras$. This result depends on a conjecture, strongly supported by results of A. Wagner  \cite{Wagner}, on the non-existence of finite non-Desarguesian planes with a simply transitive abelian group of collineations.  The relation between $\kras$-vector spaces and projective geometries also shows that, in the case of the ad\`ele class space $\H_\Q$, the hyperring structure encodes the {\em full} information on the ring structure on the ad\`eles: \cf Theorem \ref{functads} and Proposition \ref{Hopf}.\vspace{.05in}

In \cite{jamifine}, we showed that in a field endowed with a given multiplicative structure, the additive structure is encoded by a {\em bijection} $s$ of the field satisfying the two requirements that $s(0)=1$ and that $s$ commutes with its conjugates under multiplication by non-zero elements. In the same paper, we also proved that if one replaces the condition for $s$ to be a bijection by that of being a {\em retraction} (\ie $s^2 = s\circ s = s$), one obtains instead an idempotent semi-field. Therefore, it is natural to wonder if one can encode, with a similar construction, the additive structures of the hyperfield extensions of $\kras$ and $\sign$ respectively. In \S\S~ \ref{hyperequ} and \ref{hypersign} of this paper, we show that  given a multiplicative structure on a hyperfield, the additive structure is encoded by \vspace{.05in}

$(i)$~an {\em equivalence} relation commuting with its conjugates, on a hyperfield extension of $\kras$,\vspace{.05in}

$(ii)$~a partial {\em  order} relation commuting with its conjugates, on a hyperfield extension of $\sign$.\vspace{.05in}

This reformulation of the additive law in hyperfields shows that these generalized algebraic structures occupy a very natural place among the more classical notions. Along the way, we also prove  that the second axiom of projective geometry (saying that if a line meets two sides of a triangle not at their intersection then it also
    meets the third side) is equivalent to the commutativity of the equivalence relations obtained by looking at the space from  different points (\cf Lemma \ref{lemlign}). We also give an example,
using the construction of M. Hall \cite{Hall}, of an (infinite) hyperfield extension of $\kras$ whose associated geometry is a non-Desarguesian plane.

In the paper we start to investigate the content of an algebraic geometry over $\kras$. The category of commutative hyperring extensions of $\kras$ is inclusive of: algebras over fields
with semi-linear homomorphisms,  abelian groups with injective homorphisms (as explained in Proposition \ref{Lyndonlem}) and  a rather exotic land   comprising homogeneous non-Desarguesian planes.
   In \S \ref{functions}, we analyze the notion of algebraic function on  $\Spec(\H_\Q)$ defined, as in the classical case, by means of elements of the set $\Hom(\Z[T],\H_\Q)$. We use the natural coproducts $\Delta^+(T)=T\otimes 1+1\otimes T$ and $\Delta^\times(T)=T\otimes T$ on $\Z[T]$ to obtain the elementary operations on functions. \vspace{.05in}

When $\K$ is a global field, the set $P(\ads)$ of {\em prime} elements of the hyperring   $\ads$ inherits a natural structure of groupoid with the product given by multiplication and  units the set of places of $\K$. The product of two prime elements is  a prime element when the two factors sit over the same place, and over each place $v$ there exists a unique  idempotent $p_v\in P(\ads)$ (\ie $p^2_v=p_v$). The id\`ele class group $C_\K=\ads^\times$  acts by multiplication on $P(\ads)$.
When $\K$ is a function field over $\F_q$, we denote by $X$  the non singular projective algebraic curve with function field $\K$ and we let $\pi:X^{\rm ab}\to X$ be the abelian cover associated to a fixed maximal abelian extension $\K^{\rm ab}$ of $\K$. We denote by  $\Pi_1^{\rm ab}(X)$ the fundamental groupoid associated to $\pi$ and  $\Pi_1^{\rm ab}(X)'\subset \Pi_1^{\rm ab}(X)$  the subgroupoid of loops (\ie of paths whose end points coincide).  In the final part of the paper we show (Theorem \ref{ccm2prop}) that  $\Pi_1^{\rm ab}(X)'$ is
{\em canonically} isomorphic to the groupoid $P(\ads)$ and that this isomorphism is equivariant for the action of the id\`ele class group $C_\K=\ads^\times$ on $P(\ads)$ and the action of the  abelianized  Weil group on $\Pi_1^{\rm ab}(X)'$.

When ${\rm char}(\K)=0$, the above geometric interpretation is no longer available. On the other hand, the arithmetic of the hyperring $\ads$ continues to hold and the  groupoid $P(\ads)$  appears to be a natural substitute for the above groupoid of loops  and it  also supports an interpretation of the explicit formulae of Riemann-Weil.

\section{Hyperrings and hyperfields} \label{firsthyper}

In this section we shall see that the natural multiplicative monoidal structure on $\F_1=\{0,1\}$ which ignores addition can be refined, within the  category of hyperrings, to become the most basic example of a  {\em hyperfield} (\cf \cite{Krasner1}). We will refer to it as to the {\em Krasner hyperfield} $\kras$. The algebraic spectrum $\Spec\kras$ of this  hyperstructure is  the most natural lift of $\Spec\F_1$ from under $\Spec\Z$ to a basic structure mapping to $\Spec\Z$.\vspace{.05in}

In a hyperfield the additive (hyper)structure is that of a {\em canonical} hypergroup (\cf\cite{Marty} and \cite{Krasner1}).
 We start by reviewing the notion of a canonical hypergroup $H$. For our applications it will be enough to consider  this particular class of  hypergroups. We denote by $+$ the hyper-composition law in $H$. The novelty is that now the sum $x+y$ of two elements in $H$ is no longer a single element of $H$ but a non-empty subset of $H$. It is customary to  define a hyper-operation on $H$ as a map
\[
+: H\times H \to \mathcal P(H)^*
\]
 taking values into the set $\cP(H)^*$ of all non-empty subsets of $H$. Thus, $\forall a,b\in H$, $a+b$ is a non-empty subset of $H$, not necessarily a singleton. One uses the notation $\forall A,B\subseteq H,~A+B:=\{\cup (a+b)|a\in A, b\in B\}$. The definition of a  canonical hypergroup requires that $H$ has a neutral element $0\in H$ (\ie an additive identity) and that the following axioms apply:\vspace{.05in}

 $(1)$~$x+y=y+x,\qquad\forall x,y\in H$\vspace{.05in}

$(2)$~$(x+y)+z=x+(y+z),\qquad\forall x,y,z\in H$ \vspace{.05in}

$(3)$~$0+x=x= x+0,\qquad \forall x\in H$\vspace{.05in}

$(4)$~$\forall x  \in H~  \ \exists!~y(=-x)\in H\quad {\rm s.t.}\quad 0\in x+y$\vspace{.05in}

$(5)$~$x\in y+z~\Longrightarrow~ z\in x-y.$\vspace{.05in}

 The uniqueness, in $(4)$, of the symmetric element $y=-x$, for any element $x\in H$, rules out\footnote{as soon as $H$ has more than two elements}  the trivial choice of taking the addition to be the full set $H$, except for the addition with $0$.

 Property $(5)$ is usually called {\em reversibility}.
In this paper we shall always consider canonical hypergroups. \vspace{.05in}

Let $(H,+)$ be a (canonical) hypergroup and $x\in H$. The set
 $$
 O(x)=\{r\in \Z\mid \exists n\in\Z : 0\in rx+ n(x-x)\}
 $$
 is a subgroup of $\Z$.
 We say that the {\em order} of $x$ is infinite (\ie $o(x) = \infty$) if $O(x)=\{0\}$. If $o(x) \neq \infty$, the smallest positive generator  $h$ of $O(x)$ is called the {\em principal order} of $x$ (\cf \cite{Corsini} Definition 57).
Let  $q = \text{min}\{s\in\N|\exists m\neq 0, \,0\in mhx + s(x-x)\}$. The couple $(h,q)$ is then called the order of $x$.\vspace{.05in}

The notion of  a  hyperring  (\cf \cite{Krasner}, \cite{Krasner1}) is the natural generalization of the classical notion of a ring, obtained by replacing a classical addition law by a hyperaddition.

 \begin{defn}\label{hyperring}
  A hyperring  $(R,+,\cdot)$ is a non-empty set $R$ endowed with a hyperaddition $+$ and the usual multiplication $\cdot$ satisfying the following properties\vspace{.05in}

$(a)$~$(R,+)$ is a canonical hypergroup.\vspace{.05in}

$(b)$~$(R,\cdot)$ is a monoid with multiplicative identity $1$\vspace{.05in}

$(c)$~$\forall r,s,t\in R$:~~$r(s+t) = rs+rt$ and $(s+t)r = sr+tr$\vspace{.05in}

$(d)$~$\forall r\in R$:~~$r\cdot 0=0\cdot r =0$, \ie $0\in R$ is an absorbing element\vspace{.05in}

$(e)$~$0\neq 1$.

\end{defn}

In the original definition of a (Krasner) hyperring, $(R,\cdot)$ is only assumed to be a semi-group satisfying $(d)$ (\cf\cite{Davvaz1} Definition 3.1.1).\vspace{.05in}

Let $(R_1,+_1,\cdot_1)$, $(R_2,+_2,\cdot_2)$ be two hyperrings. A map $f: R_1 \to R_2$ is called a homomorphism of hyperrings if the following conditions are satisfied\vspace{.05in}

$(1)$~$f(a+_1 b)\subseteq f(a)+_2 f(b),~\forall a,b\in R_1$\vspace{.05in}

$(2)$~$f(a\cdot_1 b) = f(a)\cdot_2 f(b),~\forall a,b\in R_1.$\vspace{.05in}

The map $f$ is said to be an epimorphism if it is a surjective homomorphism such that (\cf \cite{Davvaz} Definition 2.8)
\begin{equation}\label{epi}
    x+y=\cup \{f(a+b)\mid f(a)=x,\,f(b)=y\}\qqq x,y\in R_2.
\end{equation}
It is an isomorphism if it is a bijective homomorphism satisfying $f(a+_1 b) = f(a) +_2 f(b)$, $\forall a,b\in R_1$.\vspace{.05in}

A hyperring $(R,+,\cdot)$ is called a {\em hyperfield} if $(R\setminus\{0\},\cdot)$ is a group.

\begin{defn} \label{defnkras} We denote by $\kras$ the hyperfield $(\{0,1\},+,\cdot)$ with additive neutral element $0$, satisfying the hyper-rule: $1+1=\{0,1\}$ and with the usual   multiplication, with identity $1$.

We let $\sign$ be the hyperfield $\sign=\{-1,0,1\}$ with the hyper-addition given by the ``rule of signs"
\begin{equation}\label{addsign}
    1+1=1\,, \ -1-1=-1\,, \ 1-1=-1+1=\{-1,0,1\}
\end{equation}
and the usual multiplication also given by the rule of multiplication of signs.
\end{defn}

The hyperfield $\kras$ is the natural extension, in the category of hyperrings, of the commutative (pointed) monoid $\F_1$, \ie  $(\kras,\cdot) = \F_1$. We shall refer to $\kras$ as to the Krasner hyperfield. Note that the order of the element $1\in \kras$ is the pair $(1,0)$, \ie the principal order is $1$ since $0\in 1+1$ and the secondary order is $0$ for the same reason.
In a similar manner one sees that the monoid underlying $\sign$ is $\F_{1^2}$, \ie $(\sign,\cdot) = \F_{1^2}$, where the order of the element $1\in \sign$ is the pair $(1,1)$. The homomorphism absolute value $\pi: \sign\to \kras,~\pi(x)= |x|$ is an epimorphism of hyperrings. \vspace{.05in}

 To become familiar with the operations in hyperstructures, we prove the following simple results
 \begin{prop} \label{simp} In a hyperring extension $R$ of the Krasner hyperfield $\kras$ one has
 $x+x=\{0,x\}$ for any $x\in R$ and moreover
 $$
  a\in a+b~\Longleftrightarrow~b\in\{0,a\}
 $$
 In particular, there is no hyperfield extension of $\kras$ of cardinality $3$ or $4$.
 \end{prop}

 \proof Since $1+1=\{0,1\}$ one gets $x+x=\{0,x\}$ using distributivity. Assume that $a\in a+b$ in $R$. Then since $a+a=\{0,a\}$ one has $-a=a$ so that by the reversibility condition $(6)$ in the definition of a hypergroup, one has $b\in a-a=\{0,a\}$.  Conversely, if $b\in\{0,a\}$, it follows immediately (by applying the condition $(4)$ for hypergroups) that $a\in a+b$.

 If $F$ is a hyperfield extension of $\kras$ of cardinality $3$, then $F$ contains an element $\alpha \notin \{0,1\}$. But then  one gets a contradiction since the subset $1+\alpha$ cannot contain   $0$ (since
 $1$ is its own opposite)  or $1$ or $\alpha$ (by the first part of this proposition).

 If $F$ is a hyperfield extension of $\kras$ of cardinality $4$, then let $\xi_j$ be the three non-zero elements of $F$. Then, by applying the first part of this proposition,  the sum $\xi_j+\xi_k$, for $j\neq k$ is forced to be the third non-zero element $\xi_\ell$ of $F$. This contradicts associativity of the hyperaddition for $\sum \xi_j$.
 \endproof

 Note that the above proof only uses the structure of $\kras$-vector space (\cf \S \ref{sectproj}).

\begin{rem}{\rm The same proof shows that in a hyperring extension $R$ of the   hyperfield $\sign$ one has
  \begin{equation}\label{notin}
a\in a+b~\Longleftrightarrow~b\in\{0,\pm a\}.
 \end{equation}}
\end{rem}

 Krasner gave in \cite{Krasner1} a construction of a hyperring as the quotient of a ring $R$ by a multiplicative subgroup $G$ of the group $R^\times$ of the invertible elements of $R$. This result states as follows

  \begin{prop}\label{krasner1}
  Let $R$ be a commutative ring and $G\subset R^\times$ be a subgroup of its multiplicative group. Then
  the following operations define a hyperring structure on the set $R/G$ of orbits for the action of $G$ on $R$ by multiplication\vspace{.05in}

  $\bullet$~Hyperaddition
  \begin{equation}\label{addr}
x+ y:=\left(xG+yG\right)/G \qqq x,y\in R/G
\end{equation}

$\bullet$~Multiplication
$$
xG\cdot yG=xyG\qqq x,y\in R/G.
$$
Moreover for any $x_i\in R/G$ one has
\begin{equation}\label{longsum}
    \sum x_i=\left(\sum x_i G\right)/G
\end{equation}

  \end{prop}

In particular, one can start with a field $K$ and consider the hyperring  $K/K^\times$. This way, one obtains a hyperstructure whose underlying set is made by two classes \ie the class of $0$ and that of $1$. If $K$ has more than two elements, $K/K^\times$ coincides with the Krasner hyperfield $\kras$.\vspace{.05in}

Next, we investigate in the set-up of Proposition \ref{krasner1}, under which conditions the hyperring $R/G$ contains the Krasner hyperfield $\kras$ as a sub-hyperfield.

   \begin{prop}\label{krasner2}
  Let $R$ be a commutative ring and $G\subset R^\times$ be a subgroup of the multiplicative group of units in $R$. Assume that $G\neq \{1\}$. Then,
  the hyperring $R/G$ contains $\kras$ as a sub-hyperfield if and only if $\{0\}\cup G$ is a subfield of $R$.
  \end{prop}
  \proof To verify whether $\kras\subset R/G$, it suffices to compute $1+1$ in $R/G$. By definition, $1+1$  is the union of all classes, under the multiplicative action of $G$, of elements of the form $g_1+g_2$, for $g_j\in G$ ($j=1,2$). Thus, the hyperring $R/G$ contains $\kras$ as a sub-hyperfield if and only if $G+G=\{0\}\cup G$. If this equality holds, then $\{0\}\cup G$ is stable under addition. Moreover $0\in G+G$ so that $g_1=-g_2$ for some $g_j\in G$ and thus  $-1=g_1g_2^{-1}\in G$. Thus $\{0\}\cup G$ is an additive subgroup of $R$. In fact, since $R^\times$ is a group, it follows that $G\cup\{0\}$ is a subfield of $R$. Conversely, let $F\subset R$ be a subfield and assume that $F$ is not reduced to the finite field $\F_2$. Then the multiplicative group $G=F^\times$ fulfills $G\neq \{1\}$. Moreover $G+G\subset F$ and $0\in G+G$ as $1-1=0$. Moreover, since $G$ contains at least two distinct elements $x,y$ one has $x-y\neq 0$ and thus $G+G=F$. Thus, in $R/G$ one has $1+1=\{0,1\}$ and thus $\kras\subset R/G$.
  \endproof

 \begin{example}\label{ex5}{\rm
 This simple example is an application of the above results and it shows that there exists a hyperfield extension of $\kras$ of cardinality $5$. Let $H$ be the union of $0$ with the powers of $\alpha$, $\alpha^4=1$. It is a set with 5 elements and the table of hyper-addition in $H$ is given by the following matrix
 $$
\left(
\begin{array}{ccccc}
0&1&\alpha&\alpha^2&\alpha^3\\
1 &\{0,1\} & \left\{\alpha ^2,\alpha ^3\right\} & \left\{\alpha ,\alpha ^3\right\} & \left\{\alpha ,\alpha ^2\right\} \\
\alpha& \left\{\alpha ^2,\alpha ^3\right\} & \{0,\alpha \} & \left\{1,\alpha ^3\right\} & \left\{1,\alpha ^2\right\} \\
\alpha^2& \left\{\alpha ,\alpha ^3\right\} & \left\{1,\alpha ^3\right\} & \left\{0,\alpha ^2\right\} & \{1,\alpha \} \\
 \alpha^3&\left\{\alpha ,\alpha ^2\right\} & \left\{1,\alpha ^2\right\} & \{1,\alpha \} & \left\{0,\alpha ^3\right\}
\end{array}
\right)
$$
This hyperfield structure is obtained, with $\alpha=1+\sqrt{-1}$, as the quotient of the finite field $\F_9=\F_3(\sqrt{-1})$ by the multiplicative group $\F_3^\times=\{\pm 1\}$. It follows from Proposition \ref{krasner2} that $F=\F_9/\F_3^\times$ is a hyperfield extension of $\kras$. Notice that the addition  has a very easy description since for any two distinct non-zero elements  $x,y$ the sum $x+y$ is the complement of $\{x,y,0\}$ (\cf \cite{Stratigopoulos} and Proposition \ref{Lyndonlem} below for a more general construction).
}\end{example}

The notions of ideal and prime ideal extend to the hyperring context (\cf \eg \cite{Procesi}, \cite{Davvaz1})

\begin{defn}\label{primeideal} A non-empty subset $I$ of a hyperring $R$ is called a hyperideal if\vspace{.05in}

$(a)$~$a,b\in I~\Rightarrow~a-b\subseteq I$\vspace{.05in}

$(b)$~$a\in I$, $r\in R~\Rightarrow~r\cdot a\in I.$\vspace{.05in}

The hyperideal $I\subsetneq R$ is called prime if $\forall a,b\in R$\vspace{.05in}

$(c)$~$a\cdot b\in I~\Rightarrow~a\in I$ or $b\in I$.
\end{defn}

For any hyperring $R$, we denote by $\Spec(R)$ the set of prime ideals of $R$ (\cf \cite{Procesi}). The following proposition shows that the hyperfield $\kras$ plays, among hyperrings, the equivalent role of the monoid $\F_1$ among monoids (\cf \cite{jamifine} Prop. 3.32).

\begin{prop}\label{ex1}For any hyperring $R$, the map
\[
\varphi:\Spec(R)\to \Hom(R,\kras)\,, \qquad\varphi(\ffp)=\varphi_\ffp
\]
\begin{equation}\label{phip}
    \varphi_\ffp(x)=0\qqq x\in \ffp\, , \ \ \varphi_\ffp(x)=1\qqq x\notin \ffp
\end{equation}
determines a  natural bijection of sets.
\end{prop}

\proof The map $\varphi_\ffp: R \to \kras$
is multiplicative since the complement of a prime ideal $\ffp$ in $R$ is a multiplicative set. It is compatible with the hyperaddition, using reversibility and Definition \ref{primeideal} (a). Thus the map $\varphi$
is well-defined. To define the inverse of $\varphi$, one assigns to a homomorphism of hyperrings $\rho\in \Hom (R,\kras)$  its kernel which is a prime ideal of $R$ that uniquely determines $\rho$.\endproof

Affine $\Z$-schemes, when viewed as representable functors from the category $\An$ of (commutative) rings to sets, extend canonically to the category of hyperrings as representable functors (they are represented by the same ring). This construction applies in particular to the affine line $\cD=\Sp(\Z[T])$, and one obtains $\cD(R)=\Hom(\Z[T],R)$ for any hyperring $R$. Notice though,  that while for an ordinary ring $R$, $\Hom(\Z[T],R)$ coincides with the set underlying $R$, this fact no longer holds for hyperrings. For instance,  by applying Proposition \ref{ex1} one sees that
$\cD(\kras)=\Spec(\Z[T])$ which is an infinite set unlike the set underlying $\kras$.\vspace{.05in}

To describe the elements of the set $\Hom(R,\sign)$, for any ring $R$, we first recall the definition of a symmetric cone in $R$: \cf \cite{reid}.
\begin{defn} \label{symcone} Let $R$ be a ring. A symmetric cone $P$ in  $R$ is a subset $P\subset R$ such that\vspace{.05in}

$\bullet$~ $0\notin P$, $P+P\subset P$, $PP\subset P$,\vspace{.05in}

$\bullet$~$P^c+P^c\subset P^c$ where $P^c$ is the complement of $P$ in
$R$,\vspace{.05in}

$\bullet$~$a\in P$ and $ab\in P$ imply $b\in P$,\vspace{.05in}

$\bullet$~$P-P=R$.
\end{defn}

The following proposition shows that the notion of a symmetric cone in a hyperring is equivalent to that of an element of $\Hom(R,\sign)$.

\begin{prop}\label{charhomsign}
$(1)$~A homomorphism from a  ring $R$ to the hyperring $\sign$ is determined by its kernel $\ffp\in\Spec(R)$ and a total  order on the field of fractions of the integral domain $R/\ffp$.

$(2)$~A homomorphism from a  ring $R$ to the hyperring $\sign$ is determined by a symmetric cone of $R$ in the sense of Definition \ref{symcone}.
\end{prop}

\proof $(1)$~Let $\rho\in \Hom(R,\sign)$. The kernel of $\rho$ is unchanged by composing $\rho$ with the absolute value map $\pi:\sign\to \kras$, $\pi(x) = |x|$. Thus ${\rm ker}(\rho)$ is a prime ideal $\ffp\subset R$. Moreover the map $\rho$ descends to the quotient $R/\ffp$ which is an integral domain. Let $F$ be the field of fractions of $R/\ffp$. One lets $P\subset F$ be the set of fractions of the form $x=a/b$ where $\rho(a)= \rho(b)\neq 0$. This subset of $F$ is well defined since $a/b=c/d$ means that $ad=bc$ and it follows  that $\rho(c)= \rho(d)\neq 0$. One has $\rho(0)=0$, $\rho(1)=1$ and $\rho(-1)=-1$ since $0\in \rho(1)+\rho(-1)$. Thus $P$ is also stable by addition since one can assume, in the computation of $a/b+c/d$, that $\rho(a)= \rho(b)=\rho(c)= \rho(d)=1$, so that $\rho(ad+bc)=\rho(cd)=1$. $P$ is also multiplicative. Moreover for $x\in F, x\neq 0$ one has $\pm x\in P$ for some choice of the sign. Thus $F$ is an ordered field and $\rho$ is the composition of the canonical morphism $R\to F$ with the map $F\to F/F_+^\times\sim \sign$. Conversely if one is given an order on the field of fractions of the integral domain $R/\ffp$, one can use the natural identification $F/F_+^\times\sim \sign$ to obtain the morphism $\rho$.

$(2)$~follows from $(1)$ and Theorem 2.3 of \cite{reid}. In fact one can also check directly that given a symmetric cone $P\subset R$, the following formula defines an element $\rho\in\Hom(R,\sign)$:
\begin{equation}\label{defnrhos}
    \rho(x)=\left\{
              \begin{array}{ll}
                1, & \forall x\in P\\
                -1, & \forall x\in -P \\
                0, & \hbox{otherwise.}
              \end{array}
            \right.
\end{equation}
Moreover, one easily checks that if $\rho\in\Hom(R,\sign)$ then $P=\rho^{-1}(1)$ is a symmetric cone.
\endproof
One can then apply Corollary 3.8 of \cite{reid} to obtain the following
\begin{prop}
The elements of $\cD(\sign)=\Hom(\Z[T],\sign)$ are described by
\begin{equation}\label{omega}
    \omega_\lambda(P(T))={\rm Sign}(P(\lambda))\qqq \lambda \in [-\infty,\infty]
\end{equation}
and, for $\lambda\in \bar\Q\cap \R$, by the two elements
\begin{equation}\label{omegapm}
\omega_\lambda^\pm(P(T))=\lim_{\epsilon\to 0+}{\rm Sign}(P(\lambda\pm \epsilon)).
\end{equation}
\end{prop}
\proof This follows from Corollary 3.8 of \cite{reid} for the total orders and from the first part of Proposition \ref{charhomsign} for the symmetric orders. \endproof

One can extend the above statements from the case of affine schemes to the general case (of non-affine schemes). First of all, we recall from \cite{Procesi} that to any hyperring $R$ is associated its prime spectrum $\Spec(R)$. This is a topological space endowed with a sheaf of hyperrings. Note  that it is not true for general hyperrings $R$ that the canonical map from $R$ to global sections of the structural sheaf on $\Spec R$ is bijective.

A {\em geometric hyperring space} $(X,\cO_X)$ is a topological space $X$ endowed with a sheaf of hyperrings $\cO_X$ (the structural sheaf).
As for geometric $\Z$-schemes (\cf\cite{demgab}, Chapter I, \S~1 Definition 1.1), one needs to impose the condition that the stalks of the structural sheaf of a geometric hyperring space are {\em local} algebraic structures,  \ie  they have only one maximal ideal. A homomorphism $\rho: R_1\to R_2$ of (local) hyperrings is local if the following property holds
\begin{equation}\label{locality}
\rho^{-1}(R_2^\times) = R_1^\times.
\end{equation}

A morphism $\varphi: X \to Y$ of geometric hyperring spaces is a pair $(\varphi,\varphi^\sharp)$ of a continuous map $\varphi: X\to Y$ of topological spaces and a homomorphism of sheaves of hyperrings $\varphi^\sharp: \mathcal O_Y \to \varphi_*\mathcal O_X$, 
which satisfy the property of being {\em local}, \ie  $\forall x\in X$ the homomorphisms connecting the stalks $\varphi^\sharp_x: \cO_{Y,\varphi(x)}\stackrel{}{\to}\cO_{X,x}$ are local (\cf\eqref{locality}).

With these notations we obtain the following result
\begin{prop}\label{ex2} For any  $\Z$-scheme $X$, one has a canonical identification of sets
\[
X\simeq\Hom(\Sp(\kras),X).
\]
Moreover, an element of $\Hom(\Sp(\sign),X)$ is completely determined by assigning a point $x\in X$ and a total order of the residue field $\kappa(x)$ at $x$.
\end{prop}
\proof Since $\kras$ is a hyperfield, $\{0\}\subset \kras$ is the only prime ideal and $\Spec\kras$ consists of a single point $\kappa$. Let $\rho\in \Hom(\Sp(\kras),X)$ be a morphism and $x=\rho(\kappa)\in X$. The morphism $\rho^\#$ is uniquely determined by the local morphism $\rho_x^\#:\cO_{X,x}\to \kras$. Since the ring $\cO_{X,x}$ is local, there exists only one local morphism $\rho_x^\#:\cO_{X,x}\to \kras$. Thus the map $\rho\mapsto \rho(\kappa)\in X$ is an injection from $\Hom(\Sp(\kras),X)$ to $X$. The existence of the local morphism $\cO_{X,x}\to \kras$ for any $x\in X$ shows the surjectivity. The same proof applies to describe the elements of $\Hom(\Sp(\sign),X)$ using Lemma \ref{charhomsign}.
\endproof

\section{$\kras$-vector spaces and projective geometry}\label{sectproj}

Let $R$ be a hyperring containing the Krasner hyperfield $\kras$. In this section we show, following W. Prenowitz   \cite{Prenowitz} and R. Lyndon \cite{Lyndon} that the additive hyperstructure on $R$ is entirely encoded by a projective geometry $\cP$ such that\vspace{.05in}

$\bullet$~The set of points of $\cP$ is $R\backslash 0$\vspace{.05in}

 $\bullet$~The line through two distinct points $x,y$ of $\cP$ is given by
  \begin{equation}\label{defnline}
    L(x,y)=(x+y)\cup \{x,y\}.
  \end{equation}

We shortly review the axioms of  projective geometry. They are concerned with the properties of a family
$\cL$ of subsets $L$ of a set $\cP$. The elements $L\in \cL$ are called lines. These axioms are listed as follows\vspace{.05in}

   $\P_1$: Two distinct points of $\cP$ determine a unique line $L\in\cL$ \ie
    $$
\forall x\neq y\in \cP \,, \ \exists ! \,L\in \cL\,, \ x\in L\,, \ y\in L.
    $$

    $\P_2$: If a line in $\cL$ meets two sides of a triangle not at their intersection then it also
    meets the third side, \ie
    $$
\forall x \neq y\in\cP ~\text{and}\  z\notin L(x,y),~\text{one has}$$
$$ L(y,z)\cap L(t,u)\neq \emptyset, \ \ \forall t\in L(x,y)\backslash\{x\}\,, \ u\in L(x,z)\backslash\{x\}.
    $$

   $\P_3$: Every line in $\cL$ contains at least three points.\vspace{.05in}

  We  shall consider the following small variant of the axiom $\P_3$\vspace{.05in}

  $\P'_3$: Every line in $\cL$ contains at least $4$ points.\vspace{.05in}

We  use the terminology  $\kras$-vector space to refer to a (commutative) hypergroup $E$ with a compatible action of $\kras$. Since $0\in \kras$ acts by the retraction to $\{0\}\subset E$ and $1\in\kras$ acts as the identity on $E$, the $\kras$-vector space structure on $E$  is in fact uniquely prescribed by the hypergroup structure. Thus a hypergroup $E$ is a $\kras$-vector space if and only if it fulfills the rule
\begin{equation}\label{idemcond}
   x+x=\{0,x\}\qqq x\neq 0.
\end{equation}

The next result is due essentially to W. Prenowitz   \cite{Prenowitz} and R. Lyndon \cite{Lyndon} \cf also \cite{Corsini},
 Chapter I, Theorems 30 and 34.

  \begin{prop}\label{propproj} Let $E$ be a $\kras$-vector space. Let $\cP=E\backslash\{0\}$. Then, there exists a unique geometry having $\cP$ as its set of points and satisfying \eqref{defnline}. This geometry fulfills the  three axioms $\P_1,\P_2,\P'_3$ of a projective geometry.

Conversely, let $(\cP,\cL)$ be a projective geometry fulfilling the axioms $\P_1,\P_2,\P'_3$. Let $E=\cP\cup \{0\}$ endowed with the hyperaddition having $0$ as neutral element and defined by the rule
\begin{equation}\label{defnaddi}
x+y=\left\{
      \begin{array}{ll}
       L(x,y)\backslash\{x,y\}, & \hbox{if } \ x\neq y\\
        \{0,x\}, & \hbox{if }\ x=y.
      \end{array}
    \right.
\end{equation}
Then  $E$ is a $\kras$-vector space.
  \end{prop}

  Before starting the proof of Proposition \ref{propproj} we prove the following result

  \begin{lem}\label{lemproj} Let $E$ be a $\kras$-vector space. Then for any two subsets $X,Y\subset E$ one has
\begin{equation}\label{equiinter}
    X\cap Y \neq \emptyset~ \Longleftrightarrow~ 0\in X+Y.
\end{equation}
\end{lem}

 \proof If $x\in X\cap Y$ then $0\in (x+x)\subset X+Y$. Conversely, if $0\in X+Y$, then $0\in x+y$, for some $x\in X$ and $y\in Y$. By reversibility one gets $x\in 0-y=\{y\}$ and $x=y$ so that  $X\cap Y \neq \emptyset$.\endproof

\proof(of the Proposition \ref{propproj}) We define $\cL$ as the set of subsets of $\cP=E\backslash 0$ of the form $L(x,y)=(x+y)\cup \{x,y\}$ for some $x\neq y\in \cP$. Let us check that the axiom $\P_1$ holds. We need to show that for $a\neq b\in \cP$, any line $L(x,y)$ containing $a$ and $b$ is equal to $L(a,b)$.  We show that if $z\in L(x,y)$ is distinct from $x,y$, then $L(x,z)=L(x,y)$. One has $z\in x+y$ and hence by reversibility $y\in x+z$. Thus $x+y\subset x+x+z=z\cup (x+z)$ and $L(x,y)\subset L(x,z)$. Moreover, since $y\in x+z$ one gets in the same way that $L(x,z)\subset L(x,y)$. This proves that for any two (distinct) points $a,b\in L(x,y)$ one has $L(a,b)=L(x,y)$. Indeed
  $$
  a\in L(x,y)\Rightarrow L(x,y)=L(a,x),  \ b\in L(x,y)=L(a,x)\Rightarrow L(a,b)=L(a,x)=L(x,y).
  $$
We now check the axiom $\P_2$. Let $t\in L(x,y)\backslash\{x\}\,, \ u\in L(x,z)\backslash\{x\}$. Then $x\in(y+t)\cap(u+z)$ so that by Lemma \ref{lemproj} one has $0\in y+t+u+z$. It follows again from Lemma \ref{lemproj} and the commutativity of the sum, that $(y+z)\cap(u+t)\neq \emptyset$ and $L(y,z)\cap L(t,u)\neq \emptyset$. Note that to get $x\in(y+t)\cap(u+z)$ one uses $y\neq t$ and $z\neq u$ but the validity of $\P_2$ is trivial in these cases. Thus one has $\P_2$.

By Proposition \ref{simp} one has $x\notin (x+y)$ for $x\neq y\in P$ and thus every line contains at least three points so that axiom $\P_3$ holds true. Let us show that in fact one has $\P'_3$. Assume $x+y=\{z\}$. Then $(x+y)+z=\{0,z\}$. Since $0\in x+(y+z)$ one has $x\in y+z$, but then $x\in x+(y+z)=\{0,z\}$ which is a contradiction.

Conversely, let $(\cP,\cL)$ be a projective geometry fulfilling axioms $\P_1,\P_2,\P'_3$ and  endow $E=\cP\cup \{0\}$ with the hyperaddition as in \eqref{defnaddi}. This law is associative since when $x,y,z\in \cP$ are not collinear one checks that the sum $x+y+z$ is the plane they generate with the three sides of the triangle deleted. For three distinct  points on the same line $L$, their sum is $L\cup\{0\}$ if the cardinality of the line is $>4$ and the complement of the fourth point in $L\cup\{0\}$ if the cardinality of the line is $4$.

Let us show  that $\forall x  \in E  \ \exists! y(=-x)\,, \ 0\in x+y$. We can assume $x\neq 0$. One has $0\in x+x$. Moreover for any $y\neq x$ one has $0\notin x+y=L(x,y)\backslash\{x,y\}$.

Finally we need to prove the reversibility which takes the form $x\in y+z\Rightarrow z\in x+y$. If $y=0$ or $z=0$, the conclusion is obvious, thus we can assume that $y,z\neq 0$. If $y=z$ then $y+z=\{0,z\}$ and one gets $z\in x+y$. Thus we can assume $y\neq z$. Then $x\in y+z$ means that $x\in L(y,z)\backslash\{y,z\}$
and this implies $z\in L(x,y)\backslash\{x,y\}$.
\endproof

\begin{rem}\label{finitedim} {\rm  Let $V$ be a $\kras$-vector space. For any finite subset $F=\{x_j\}_{j\in J}\subset \cP=V\setminus 0$, the subset $$E=\{\sum_{j\in J}\lambda_j x_j\mid \lambda_j\in \kras\}$$  is stable under hyperaddition and it follows from the formula $(x+x)=\{0,x\}$ that $E$ coincides with $\displaystyle{\sum_{j\in J}}(x_j+x_j)$.   Thus, $W=E\setminus 0$ is a subspace of the geometry $\cP$  \ie  a subset of $V\setminus 0$ such that
\begin{equation}\label{subsp}
   \forall x\neq y\in  W\, \ \ L(x,y)\subset W
\end{equation}
and the restriction to $W=E\setminus 0$ of the geometry of $\cP$ is finite dimensional. We refer to \cite{Stratigopoulos} for  the notion  of dimension of a vector space over a hyperfield. Here, such dimension is related to the dimension $\dim W$ of the associated projective geometry by the equation
\begin{equation}\label{subsp1}
   \dim W=\dim_\kras E -1\,.
\end{equation}
}\end{rem}

Next result shows that hyperfield extensions of $\kras$  correspond precisely to the ``Zweiseitiger  Inzidenszgruppen" (two-sided incidence groups) of \cite{Ellers}. In particular, the {\em commutative} hyperfield extensions of $\kras$  are classified by  projective geometries together with a simply transitive action by a commutative subgroup of the collineation group. We first recall the definition of a two-sided incidence group

\begin{defn}\label{Inzidenszgruppen} Let $G$ be a group which is the set of points of a projective geometry.
Then $G$ is called a two-sided incidence group if the left and right translations by $G$ are automorphisms of the geometry.
\end{defn}

We can now state the precise relation between  hyperfield extensions of $\kras$ and two-sided incidence groups (\cite{Ellers} and \cite{Ellers1}) whose projective geometry satisfies the axiom $\P'_3$ in place of  $\P_3$.

\begin{lem}\label{main} Let $\H\supset \kras$ be a  hyperfield extension of $\kras$. Let $(\cP,\cL)$ be the associated geometry (\cf Proposition  \ref{propproj}). Then, the multiplicative group $\H^\times$ endowed with the geometry $(\cP,\cL)$ is a two-sided incidence group fulfilling $\P'_3$.

Conversely, let $G$ be a two-sided incidence group fulfilling $\P'_3$. Then, there exists a unique hyperfield extension $\H\supset \kras$ such that $\H=G\cup \{0\}$. The hyperaddition in $\H$ is defined by the rule $$x+y=L(x,y)\backslash\{x,y\}\quad \text{for any}~x\neq y\in \cP$$ and the multiplication is the group law of $G$, extended by $0\cdot g=g\cdot 0=0$,  $\forall g\in G$.
\end{lem}

\proof For the proof of the first statement it suffices to check that the left and right multiplication by a non-zero element $z\in \H$ is a collineation. This follows from the distributivity property of the hyperaddition which implies that
 \begin{equation}\label{added-formula}
zL(x,y)=z(x+y)\cup \{zx,zy\}=L(zx,zy)\,.
\end{equation}
A similar argument shows that the right multiplication is also a collineation.

Conversely, let $G$ be a two-sided incidence group fulfilling $\P'_3$.  Let $\H=G\cup \{0\}$ and define the hyperaddition as in Proposition \ref{propproj}.
With this operation, $\H$ is an additive hypergroup. Let the multiplication be the group law of $G$, extended by $0\cdot g=g\cdot 0=0$, $\forall g\in G$.
This operation is distributive with respect to the hyperaddition because $G$ acts by collineations. Thus one obtains an hyperfield $\H$. Moreover, by construction, the projective geometry underlying $\H$ is $(\cP,\cL)$. \endproof

 Let $H$ be an abelian group. We define the geometry on $H$ to be that of a single line.  By applying Lemma \ref{main}, we obtain the following result (\cf \cite{Stratigopoulos}, Proposition 2)
 \begin{prop} \label{Lyndonlem} Let $H$ be an abelian group of order at least $4$. Then, there exists a unique hyperfield extension
$\kras[H]$ of $\kras$ whose underlying monoid is $\F_1[H]$  and whose geometry is that of a single line.

The assignment $H\mapsto \kras[H]$ is  functorial only for  injective homomorphisms of abelian groups and for the canonical surjection $\kras[H]\to\kras$.
 \end{prop}
\proof Let $R=H\cup\{0\}$ viewed as a monoid. The construction of Lemma \ref{main} gives the following hyperaddition on $R=\kras[H]$ (\cf \cite{Lyndon})
\begin{equation}\label{addlyndon}
    x+y=\left\{
         \begin{array}{ll}
           x, & \hbox{if}\  y=0\\
           \{0,x\}, & \hbox{if}\ y=x \\
           R\backslash\{0,x,y\}, & \hbox{if}\ \#\{0,x,y\}=3.
         \end{array}
       \right.
\end{equation}
One easily checks that this (hyper)operation determines a hypergroup law on $\kras[H]$, provided that the order of $H$ is at least $4$. Moreover, since the left multiplication by a non-zero element is a bijection preserving $0$, one gets the distributivity.  Let then
$\rho:H_1\to H_2$ be a group homomorphism. If $\rho$ is injective and $x\neq y$ are elements of $H_1$ then,
by extending $\rho$ by $\rho(0)=0$, one sees that $\rho(x+y)\subset\rho(x)+\rho(y)$.  If $\rho$ is not injective and does not factor through $\kras[H_1]\to\kras\subset \kras[H_2]$, then there exists two elements of $H_1$, $x\neq y$ such that $\rho(x)=\rho(y)\neq 1$. This contradicts the required property  $\rho(x+y)\subset\rho(x)+\rho(y)$ of a homomorphism of hyperrings (\cf\S~\ref{firsthyper}) since $\rho(x)+\rho(y)=\{0,\rho(x)\}$ while $1\in x+y$ so that $1=\rho(1)\in\rho(x+y)$.\endproof

\begin{rem}\label{monoids}{\rm The association $H\to \kras[H]$ determines a functor from abelian groups (and injective morphisms) to hyperfield extensions of $\kras$. This functor does not extend to  a functor from monoids to hyperring extensions of $\kras$  since the distributivity (of left/right multiplication) with respect to the addition \eqref{added-formula} fails in general when $H$ is only a monoid.

One can show that all commutative hyperring extensions $R$ of $\kras$ such that $\dim_\kras(R)=2$ are of the form $R=\kras[H]^{(j)}$ for some $j\in \{0,1,2\}$ where $H$ is an abelian group of cardinality $>3-j$. Here $\kras[H]^{(0)}=\kras[H]$, $\kras[H]^{(1)}=\kras[H]\cup \{a\}$ with the presentation
\begin{equation}\label{(1)}
     a^2=0,\ au=ua=a\qqq u\in H
\end{equation}
while $ \kras[H]^{(2)}=\kras[H]\cup \{e,f\}$ with the presentation (\cf \cite{Stratigopoulos})
\begin{equation}\label{(2)}
    e^2=e,\ f^2=f, \ ef=fe=0,\ au=ua=a\qqq u\in H, a\in\{e,f\}\,.
\end{equation}
}\end{rem}

Next result is, in view of Lemma \ref{main}, a  restatement of the classification of Desarguesian ``Kommutative Inzidenszgruppen" of \cite{Karzel}.

\begin{thm}\label{thmmain} Let $\H\supset \kras$ be a commutative hyperfield extension of $\kras$. Assume that the geometry associated to the $\kras$-vector space $\H$ is Desarguesian\footnote{This is automatic if the $\dim_\kras\H$ is $>3$}
and of dimension $\geq 2$. Then, there exists a  unique pair $(F,K)$ of a commutative field $F$ and a subfield $K\subset F$ such that
\begin{equation}\label{isomain}
    \H=F/K^\times.
\end{equation}
\end{thm}

\proof By applying Lemma \ref{main} one gets a  Desarguesian geometry with a simply transitive action of an abelian group by collineations. It follows from \cite{Karzel} (\S 5 Satz 3) that there exists a normal near-field $(F,K)$ such that the commutative incidence group is  $F^\times/K^\times$. By \opcit (\S 7 Satz 7), the near-field $F$ is in fact a commutative field. The uniqueness of this construction follows from \opcit: \S 5, (5.8).\endproof

\begin{rem}\label{ncremark}{\rm By applying the results of H. W\"{a}hling (\cf\cite{Wahling}), the above Theorem \ref{thmmain} is still valid without the hypothesis of commutativity (for the multiplication) of $\H$. The field $F$ is then a skew field and $K$ is {\em central}  in $F$.
}\end{rem}

Theorem~\ref{thmmain} generalizes  to the case of commutative, integral hyperring extensions of $\kras$.

\begin{cor}\label{corthmmain} Let $\H\supset \kras$ be a commutative hyperring extension of $\kras$. Assume that $\H$ has no zero divisors and that $\dim_\kras\H >3$. Then, there exists a unique pair $(A,K)$ of a commutative integral domain $A$ and a subfield $K\subset A$ such that
\begin{equation}\label{isomain1}
    \H=A/K^\times.
\end{equation}
\end{cor}

\proof By \cite{Procesi}, Prop. 6 and 7 (\cf also \cite{Davvaz}), $\H$ embeds in its hyperfield of fractions. Thus, by applying Theorem \ref{thmmain} one obtains the desired result.\endproof

\subsection{Finite extensions of $\kras$}

In view of Theorem \ref{thmmain},  the classification of all finite, {\em  commutative}  hyperfield extensions of $\kras$ reduces to the determination of non-Desarguesian finite projective planes with a simply transitive abelian group $G$ of collineations.

\begin{thm}\label{classcor}
Let $\H\supset \kras$ be a finite commutative  hyperfield extension of $\kras$. Then,  one of the following cases occurs\vspace{.05in}

$(1)$~$\H=\kras[G]$ (\cf Proposition \ref{Lyndonlem}), for a finite abelian group $G$.\vspace{.05in}

$(2)$~There exists a finite field extension $\F_q\subset \F_{q^m}$ of a finite field $\F_q$ such that $\H=\F_{q^m}/\F_q^\times$.\vspace{.05in}

$(3)$~There exists a finite, non-Desarguesian projective plane $\cP$ and a simply transitive abelian group $G$ of collineations of $\cP$, such that $G$ is the commutative incidence group associated to $\H$ by Lemma  \ref{main}.
 \end{thm}

\proof Let $G$ be the incidence group associated to $\H$ by Lemma  \ref{main}. Then, if the geometry on $G$ consists of a single line, case (1) applies. If the geometry associated to $\H$ is Desarguesian, then by Theorem \ref{thmmain}  case (2) applies. If neither $(1)$ nor $(2)$ apply, then the geometry of $\H$ is a finite non-Desarguesian plane with a simply transitive abelian group $G$ of collineations. \endproof

\begin{rem}{\rm There are no known examples of finite, commutative hyperfield extensions $\H\supset\kras$ producing  projective planes as in case (3) of the above theorem. In fact, there is a conjecture (\cf \cite{Beutel} page 114) based on results of  A. Wagner and T. Ostrom (\cf \cite{Beutel} Theorem 2.1.1, Theorem 2.1.2, \cite{Wagner}, \cite{Wagner1}) stating that such case cannot occur. A recent result of K. Thas and D. Zagier \cite{Thas} relates  the existence of potential counter-examples to  Fermat curves and surfaces and number-theoretic exponential sums. More precisely, the existence of a special prime $p=n^2+n+1$ in the sense of \opcit Theorem 3.1 (other than $7$ and $73$) would imply the existence of a non-Desarguesian plane $\Pi=\Pi(\F_p,(\F_p^\times)^n)$ as in case $(3)$ of the above theorem. Note that, by a result of M. Hall \cite{Hall} there exists {\em infinite} non-Desarguesian projective planes with a {\em cyclic} simply transitive group of collineations. We shall come back to the corresponding hyperfield extensions of $\kras$ in \S \ref{hyperequ}. }
\end{rem}

\subsection{Morphisms of quotient hyperrings}
\label{endomorphisms}

 Let $E,F$ be   $\kras$-vector spaces. Let $T:E\to F$ be a homomorphism of hypergroups (respecting the action of $\kras$).  The kernel of $T$
$$
\Ker\, T=\{\xi\in E\mid T\xi=0\}
$$
intersects $\cP_E=E\backslash \{0\}$ as a subspace $N=\Ker\, T\cap\cP_E$ of the geometry  $(\cP_E,\cL_E)$. For any $\eta\in \cP_E$, the value of $T(\eta)$ only depends  upon the subspace $N(\eta)$ of $\cP_E$ generated by $N$ and $\eta$, since $T(\eta+\xi)\subset T(\eta)+T(\xi)=T(\eta)$ for $\xi\in N$. One obtains in this way a morphism of projective geometries in the sense of \cite{Faure} from $(\cP_E,\cL_E)$  to $(\cP_F,\cL_F)$. More precisely the restriction of $T$ to the complement of $\Ker\, T$ in $\cP_E$ satisfies the following properties\vspace{.05in}

$(M1)$~$N$ is a subspace of $\cP_E$.

$(M2)$~$a,b\notin N$, $c\in N$ and $a\in L(b,c)$ imply $T(a)=T(b)$.\vspace{.05in}

$(M3)$~$a,b,c\notin N$ and $a\in b\vee c$ imply $T(a)\in T(b)\vee T(c)$.\vspace{.05in}

In the last property one sets $x\vee y=L(x,y)$ if $x\neq y$ and $x\vee y=x$ if $x=y$. Note that $(M3)$ implies that if $T(b)\neq T(c)$ the map $T$ injects the line $L(b,c)$ in the line $L(T(b),T(c))$.

Conversely one checks that any morphism of projective geometries (fulfilling $\P_3'$) in the sense of \cite{Faure} comes from a unique morphism of the associated $\kras$-vector spaces.

  A complete description of the non-degenerate\footnote{a morphism is non-degenerate when its range is not contained in a line}  morphisms of Desarguesian geometries in terms of semi-linear maps is also given in   \opcit  In our context we use it to show the following result

\begin{thm}\label{main0} Let $A_j$ ($j=1,2$) be a commutative algebra over the field $K_j\neq \F_2$, and let
$$
\rho\,:\, A_1/K_1^\times\to A_2/K_2^\times
$$
be a homomorphism of hyperrings. Assume that the range of $\rho$ is of $\kras$-dimension $>2$, then $\rho$ is induced by a unique ring homomorphism $\tilde \rho:A_1\to A_2$ such that $\alpha=\tilde \rho|_{K_1}$ is a field inclusion $\alpha:K_1\to K_2$.
\end{thm}
\proof  Since $\rho$ is a homomorphism of $\kras$-vector spaces, it defines a morphism of projective geometries in the sense of \cite{Faure}. Moreover,  since $\rho$ is non-degenerate by hypothesis, there exists by \opcit Theorem 5.4.1, (\cf also \cite{Faure1} Theorem 3.1), a semi-linear map $f:A_1\to A_2$ inducing $\rho$. We let $\alpha:K_1\to K_2$ be the corresponding morphism of fields. Moreover, $f$ is uniquely determined up to multiplication by a scalar, and hence it is uniquely fixed by the condition $f(1)=1$ (which is possible since $\rho(1)=1$ by hypothesis). Let us show that, with this normalization, the map $f$ is a homomorphism. First of all, since $\rho$ is a homomorphism one has
\begin{equation}\label{prop}
f(xy)\in K_2^\times f(x)f(y)\qqq x,y\in A_1\,.
\end{equation}
Let us then show that if $\rho(x)\neq 1$ one has $f(xy)=f(x)f(y)$ for all $y\in A_1$. We can assume, using \eqref{prop}, that $f(x)f(y)\neq 0$ and we let $\lambda_{x,y}\in K_2^\times$ be such that $f(xy)=\lambda_{x,y} f(x)f(y)$. We assume $\lambda_{x,y}\neq 1$ and get a contradiction. Let us show that
\begin{equation}\label{prop1}
    \xi(s,t)=1+\alpha(s)f(x)+\alpha(t)f(y)\in K_2f(x)f(y)\qqq s,t\in K_1^\times.
\end{equation}
This follows from \eqref{prop} which proves the collinearity of the vectors
$$
f\left((1+sx)(1+ty)\right)=\xi(s,t)+\alpha(st)f(xy)=\xi(s,t)+\alpha(st)\lambda_{x,y} f(x)f(y)
$$
$$
f(1+sx)f(1+ty)=\xi(s,t)+\alpha(st)f(x)f(y)
$$
Thus by \eqref{prop1} the vectors $\xi(s,t)$ are all proportional to a fixed vector. Taking two distinct $t\in K_1^\times$ shows that $f(y)$ is in the linear span of the (independent) vectors $1,f(x)$ \ie $f(y)=a+bf(x)$ for some $a,b\in K_2$. But then taking $t$ with $1+\alpha(t)a \neq 0$ and two distinct $s\in K_1^\times$ contradicts the proportionality since $1$, $f(x)$ are independent, while
$$
\xi(s,t)=(1+\alpha(t)a)1+(\alpha(s)+\alpha(t)b)f(x)\,.
$$
Thus we have shown that if $\rho(x)\neq 1$ one has $f(xy)=f(x)f(y)$ for all $y\in A_1$.
Let then $x_0\in A_1$ be such that $\rho(x_0)\neq 1$. One has $f(x_0y)=f(x_0)f(y)$ for all $y\in A_1$. Then for $x\in A_1$ with $\rho(x)=1$ one has $\rho(x+x_0)\neq 1$ and the equality $f((x+x_0)y)=f(x+x_0)f(y)$ for all $y\in A_1$ gives $f(xy)=f(x)f(y)$.
\endproof

\begin{cor}\label{main1} Let $A$ and $B$ be commutative algebras over $\Q$ and let
$$
\rho\,:\, A/\Q^\times\to B/\Q^\times
$$
be a homomorphism of hyperrings. Assume that the range of $\rho$ is of $\kras$-dimension $>2$, then $\rho$ is induced by a unique ring homomorphism $\tilde \rho:A\to B$.
\end{cor}

\begin{rem}\label{quot}{\rm Let $A$ and $B$ be commutative  $\Q$-algebras and let
$$
\rho\,:\, A \to B/\Q^\times
$$
be a homomorphism of hyperrings. One has $\rho(1)=1$ by hypothesis. By induction one gets $\rho(n)\in \{0,1\}$
for $n\in \N$. Moreover, since $0=\rho(0)\in \rho(1)+\rho(-1)$,  $\rho(-1)$ is the additive inverse of $1$ in $B/\Q^\times$,  it follows that  $\rho(-1)=1$. By the multiplicativity of $\rho$ one gets $\rho(n)\in \{0,1\}$
for $n\in \Z$. Using the property that $n\cdot 1/n=1$ one obtains $\rho(n)=1$
for $n\in \Z, n\neq 0$. Again  by the multiplicativity of $\rho$, it follows that $\rho$ induces a homomorphism
$A/\Q^\times \to B/\Q^\times$ and Corollary \ref{main1} applies.
}\end{rem}

\begin{rem}\label{sub2}{\rm Let $A,B,\rho$ be as in Corollary \ref{main1}.
 Assume that the range $\rho(A)$ of $\rho$ has $\kras$-dimension $\leq 2$.
  Then, one has $\rho(1)=1\in \rho(A)$, and either $\rho(A)=\kras$ or
  there exists $\xi \in \rho(A)$, $\xi\notin \kras$ such that $\rho(A)\subset \Q+\Q b\subset B/\Q^\times$
  where $b\in B$ is a lift of $\xi$. Since $\rho$ is multiplicative, one has
  $\xi^2\in \rho(A)$ and
   $b$ fulfills a quadratic equation
$$
b^2=\alpha+\beta b\,, \ \alpha,\beta\in \Q.
$$
One can reduce to the case when $b$ fulfills the condition
\begin{equation}\label{quadequ}
   b^2=N\,, \ N\in \Z\,, \ N\  \text{square\, free}.
\end{equation}
Thus the morphism $\rho\,:\, A/\Q^\times\to B/\Q^\times$ factorizes through the
quadratic subalgebra $\Q(\sqrt N):=\Q[T]/(T^2-N)$
\begin{equation}\label{factorrho}
   \rho\,:\, A/\Q^\times\to \Q(\sqrt N)/\Q^\times \to B/\Q^\times.
\end{equation}
Let us consider the case $N=1$. In this case $\Q(\sqrt N)$ is  the algebra $B_0=\Q\oplus \Q$  direct sum of two copies of $\Q$. For $n\in \N$, an odd number, the map $\rho_n: B_0\to B_0,~\rho_n(x)=x^n$ is a multiplicative endomorphism of $B_0$. Let $\tilde \P_1=B_0/\Q^\times$ be the quotient hyperring. The corresponding geometry is the projective line $\P^1(\Q)$ and for any $x\neq y\in \tilde \P_1\backslash \{0\}$ one has
$$
x+y=\tilde \P_1\backslash \{0,x,y\}.
$$
Since $\rho_n$ induces an injective  self-map of $\tilde \P_1$, one gets that
$$
\rho_n(x+y)=\rho_n(\tilde \P_1\backslash \{0,x,y\})\subset \tilde \P_1\backslash \{0,\rho_n(x),\rho_n(y)\}
=\rho_n(x)+\rho_n(y).
$$
Thus $\rho_n: B_0/\Q^\times\to B_0/\Q^\times$ is an example of a morphism of hyperrings which does not lift to a ring homomorphism. The same construction applies when the map  $x\mapsto x^n$ is replaced by any injective group homomorphism $\Q^\times\to \Q^\times$.
}\end{rem}

\section{The equivalence relation on  a hyperfield extension of $\kras$} \label{hyperequ}

In this section we prove that the addition in a hyperfield extension $F$ of the Krasner hyperfield $\kras$ is uniquely determined by  an equivalence relation on $F$ whose main property is that to commute with its conjugates by rotations.

\subsection{Commuting  relations}

Given two  relations $T_j$ $(j=1,2$) on a set $X$, one defines  their composition as
$$
T_1\circ T_2=\{(x,z)\mid \exists y\in X, \ (x,y)\in T_1\,, \ (y,z)\in T_2\}.
$$
By definition, an equivalence relation $T$ on a set $X$ fulfills $\Delta\subset T$, where $\Delta=\Delta_X$ denotes the diagonal. Moreover, one has $T^{-1}=T$ where
$$
T^{-1}=\{(x,y)\mid (y,x)\in T\}
$$
and finally $T\circ T=T$.

We say that two equivalence relations $T_j$ on a set $X$ commute when any of the following equivalent conditions hold:\vspace{.05in}

$\bullet$~$T_1\circ T_2=T_2\circ T_1$\vspace{.05in}

$\bullet$~$T_1\circ T_2$ is the equivalence relation generated by the $T_j$\vspace{.05in}

$\bullet$~$T_1\circ T_2$ is an equivalence relation.\vspace{.05in}

 Notice that any of the above conditions  holds if and only if for any class $C$ of the equivalence relation generated by the $T_j$,
the restrictions of the $T_j$ to $C$ are independent in the sense that any class of $T_1|_C$ meets every class of $T_2|_C$.

\subsection{Projective geometry as commuting points of view}

Given a point $a\in \cP$ in a projective geometry $(\cP,\cL)$, one gets a natural equivalence relation $R_a$ which partitions the set of points $\cP\backslash\{a\}$  as the lines through $a$. We extend this to an equivalence relation  $R_a$, denoted $\sim_a$, on $\cP\cup\{0\}$ such that $0\sim_a a$ and for $x\neq y$ not in $\{0,a\}$
    \begin{equation}\label{aligned}
        x\sim_a y~\Longleftrightarrow~ a\in L(x,y).
    \end{equation}
    We now relate the commutativity of these equivalence relations with the axiom $\P_2$. More precisely, we have the following

    \begin{lem}\label{lemlign} The axiom $\P_2$ of a projective geometry $(\cP,\cL)$ is equivalent to the commutativity of the equivalence relations $R_a$.
    \end{lem}
    \proof Let us first assume that the axiom  $\P_2$ holds and show that the relations $R_a$'s commute pairwise. Given two points $a\neq b$ in $\cP$, we first determine the equivalence relation $R_{ab}$ generated by $R_a$ and $R_b$. We claim that the equivalence classes for $R_{ab}$ are\vspace{.05in}

$\bullet$~The union of $L(a,b)$ with $\{0\}$.\vspace{.05in}

$\bullet$~The complement of $L(a,b)$ in any plane containing $L(a,b)$.\vspace{.05in}

          One checks indeed that these subsets are stable under  $R_a$ and $R_b$. Moreover let us show that in each of these subsets, an equivalence class of $R_a$ meets each equivalence class of $R_b$. In the first case, $R_a$ has two classes: $\{0,a\}$ and $L(a,b)\backslash \{a\}$  (similarly for $R_b$), so the result is clear. For the complement of $L(a,b)$ in any plane containing $L(a,b)$, each class of $R_a$ is the complement of $a$ in a line through $a$ and thus meets each class of $R_b$, since coplanar lines meet non-trivially. Thus $R_a$ commutes with $R_b$.

     Conversely, assume that for all $a\neq b$ the relation $R_a$ commutes with $R_b$.   Let then $x,y,z,t,u$ as in the statement of the axiom $\P_2$. One has
     $t\sim_y x$ and $z\sim_u x$. Thus $z\in R_uR_y(t)$.  Then $z\in R_yR_u(t)$  and
     $ L(y,z)\cap L(t,u)\neq \emptyset$.\endproof

     We can thus reformulate the axioms of projective geometry in terms of a collection of {\em commuting points of view}, more precisely:

    \begin{prop} \label{propequigeom} Let $X=\cP\cup \{0\}$ be a pointed set and let $\{R_a; a\in \cP\}$ be a family of equivalence relations on $X$ such that\vspace{.05in}

$(1)$~$R_a$ commutes with $R_b$,  $\forall~a,b\in \cP$\vspace{.05in}

$(2)$~$\{0,a\}$ is an equivalence class for $R_a$, for all $a\in \cP$\vspace{.05in}

$(3)$~Each equivalence class of $R_a$, other than $\{0,a\}$, contains at least three elements.\vspace{.05in}

For $a\neq b\in \cP$ let $L(a,b)$ be the intersection with $\cP$ of the  class of $0$ for $R_a\circ R_b$.
Define a collection $\cL$ of lines in $\cP$ as the set of all lines $L(a,b)$. Then  $(\cP,\cL)$ is a projective geometry fulfilling the axioms $\P_1$, $\P_2$ and $\P'_3$.
\end{prop}

 \proof One has $R_b(0)=\{0,b\}$ and thus the points of $L(a,b)\backslash \{a\}$ are those of $R_a(b)$. The same statement holds after interchanging $a$ and $b$.
Let us  show that if $c\in L(a,b)$ is distinct from both $a$ and $b$, then $L(a,c)=L(a,b)$. The points of $L(a,c)\backslash \{a\}$ are those of $R_a(c)$ and $c\in R_a(b)$ since $c\in L(a,b)\backslash \{a\}$. By transitivity it follows $R_a(c)=R_a(b)$. Thus $L(a,c)\backslash \{a\}=L(a,b)\backslash \{a\}$ and $L(a,c)=L(a,b)$.
Hence, for any two (distinct) points $a,b\in L(x,y)$ one has $L(a,b)=L(x,y)$. Thus, if we let the set $\cL$ of lines in $\cP$ be given by all $L(a,b)$ axiom $\P_1$ follows while the condition (3) ensures $\P'_3$. For $x\neq y$ not in $\{0,a\}$,  one has that $a\in L(x,y)$ iff $x\in R_a(y)$. Indeed, if $a\in L(x,y)$ then $L(x,y)=L(y,a)$ and $x\in L(y,a)\backslash \{a\}=R_a(y)$. Conversely, if $x\in R_a(y)$, then $x\in L(y,a)$ and $a\in L(x,y)$. Thus, by Lemma \ref{lemlign} one gets $\P_2$. \endproof

\subsection{The basic equivalence relation on  a hyperfield extension of $\kras$} \label{subhyperequ}

In the case of a hyperring containing $\kras$, the following statement shows that the equivalence relation associated to the unit $1$ plays a privileged role.

\begin{prop}\label{addequ} Let $R$ be a hyperring containing $\kras$ as a sub-hyperring. We introduce the multi-valued map $s: R \to R$, $s(a) = a+1$. Then, the following conditions are equivalent. For $x,y\in R$\vspace{.05in}

$(1)$~$x=y$ or $x\in y+1$\vspace{.05in}

$(2)$~$x\cup (x+1)=y\cup (y+1)$\vspace{.05in}

$(3)$~$s^2(x)=s^2(y)$, ($s^2 = s\circ s$).\vspace{.05in}

The above equivalent conditions define an equivalence relation $\sim_R$ on $R$.
\end{prop}

\proof We show that (1) implies (2). Assume $x\in y+1$. Then $x+1\subset y+1+1=y\cup (y+1)$. Thus
$x\cup (x+1)\subset y\cup (y+1)$. By reversibility one has $y\in x+1$ and thus
$y\cup (y+1)\subset x\cup (x+1)$ so that $x\cup (x+1)= y\cup (y+1)$.

Next, we claim that (2) and (3) are equivalent since $s^2(a)=a+1+1=a\cup (a+1)$ for any $a$. Finally (2) implies (1), since if $x\neq y$ and (2) holds one has $x\in y+1$.\endproof

 One knows by Proposition \ref{simp} that $a\notin s(a)$ provided that $a\neq 1$. It follows that the map $s$ is in fact completely determined by the equivalence relation $\sim_R$. Thus one obtains

\begin{cor}\label{addequ1} Let $R$ be a hyperring containing the Krasner hyperfield $\kras$ and let $\sim_R$ be the associated equivalence relation.  Then one has
\begin{equation}\label{sxequ}
    x+1=\{y\sim_R x\,, \ y\neq x\}\qqq x\in R,\ x\neq 1.
\end{equation}
 In particular, when $R$ is a hyperfield its additive hyper-structure is uniquely determined by the equivalence relation $\sim_R$.
\end{cor}

We now check directly the commutativity of $\sim_R$ with  its conjugates under multiplication by any element $a\in R^\times$.

\begin{lem}\label{comrel} Let $R$ be a hyperring containing the Krasner hyperfield $\kras$ as sub-hyperring and let $\sim_R$ be the corresponding equivalence relation.  Then $\sim_R$ commutes with its conjugates under multiplication by any element $a\in R^\times$.
\end{lem}
\proof Let $T=\sim_R$. One has $T(x)=x+1+1$ for all $x\in R$. It follows that for the conjugate relation $T^a:=aTa^{-1}$ one has $T^a(x)=x+a+a$. Thus
$$
T\circ T^a(x)=1+1+(a+a+x)=a+a+(1+1+x)=T^a\circ T(x).
$$
\endproof
Thus, one can start with any abelian group $H$ (denoted multiplicatively) and by applying  Corollary \ref{addequ1}, consider on the set $R=H\cup \{0\}$ an equivalence relation $S$ which commutes with its conjugates under rotations. Let us assume that $\{0,1\}$ forms an equivalence class for $S$. In this generality, it is not true that   the multivalued map $s: R \to R$ defined by
\begin{equation}\label{sxequbis}
    s(x)=\{y\in S(x)\,, \ y\neq x\}\qqq x\in R,\ x\neq 1,\qquad s(1)=\{0,1\}
\end{equation}
commutes with its conjugates under rotations. One can consider, for example, $H=\Z/3\Z$ and on the set $R=H\cup \{0\}$ one can define the equivalence relation $S$ with classes $\{0,1\}$ and $\{j,j^2\}$. This relation $S$ commutes with its conjugates under rotations, but one has
$$
s^j(s(1)=\{j,j^2\}\,, \ \ s(s^j(1))=j.
$$
But the commutativity of $s$ with its conjugates holds provided the equivalence classes for $S$ other than $\{0,1\}$ have cardinality at least three.  One in fact obtains the following
\begin{prop}\label{lemadd} Let  $H$ be an abelian group.
Let $S$ be an equivalence relation on the set $R=H\cup \{0\}$ such that\vspace{.05in}

$\bullet$~$\{0,1\}$ forms an equivalence class for $S$\vspace{.05in}

$\bullet$~Each class of $S$, except $\{0,1\}$, contains at least three elements\vspace{.05in}

$\bullet$~The relation $S$ commutes with its conjugates   for the action of $H$ by multiplication on the monoid $R$.\vspace{.05in}

Then with $s$ defined by \eqref{sxequbis}, the law
\begin{equation}\label{mainadddefn}
    x+y:=\begin{cases}  y&\text{if}~x=0\\
xs(yx^{-1})&\text{if}~x\neq 0\end{cases}
\end{equation}
defines a commutative hypergroup structure on $R$. With this hyper-addition the monoid   $R$ becomes a commutative hyperfield containing $\kras$.
\end{prop}

\proof For each $a\in H$ let $S_a$ be the equivalence relation obtained by conjugating $S$ by the multiplication by $a$. Thus $x\sim y\,, \ (S_a)$ means $a^{-1}x\sim a^{-1}y\,, \ (S)$. In particular $\{0,a\}$ is an equivalence class for $S_a$. One checks that all conditions of Proposition \ref{propequigeom} are fulfilled and thus one gets a geometry fulfilling axioms $\P_1$, $\P_2$ and $\P'_3$.
By construction the abelian group $H$ acts by collineations on this geometry and thus Theorem \ref{main} applies. \endproof

Note that one can give a direct proof of Proposition \ref{lemadd}, in fact we shall use that approach to treat a similar case in \S \ref{hypersign}.

\begin{example}{\rm
The construction of projective planes from {\em difference sets} (\cf \cite{Singer}) is a special case of Proposition \ref{lemadd}. Let $H$ be an abelian group, and $\cD\subset H$ be a subset of $H$ such that the following map is bijective
$$
\cD\times \cD\backslash \Delta\to H\backslash \{1\},\ (x,y)\mapsto xy^{-1}
$$
(where $\Delta$ is the diagonal). Then the partition of $H\backslash \{1\}$ into the subsets $\cD\times \{u\}$ for $u\in \cD$ defines on $R=H\cup \{0\}$ an equivalence relation $S$ which fulfills all conditions of Proposition \ref{lemadd}. By \cite{Hall} Theorem 2.1 one obtains in this manner all cyclic projective planes \ie in the above context all hyperfield extensions of $\kras$ whose multiplicative group is cyclic and whose associated geometry is of dimension $2$. By \cite{Hall} Theorem 3.1, difference sets $\cD$ exist for the infinite cyclic group $\Z$ and thus provide examples of hyperfield extensions of $\kras$ whose multiplicative group is cyclic and whose associated geometry is non-Desarguesian.
}\end{example}

\section{The  order relation on  a hyperfield extension of $\sign$} \label{hypersign}\vspace{.05in}

Let $\sign$ be the hyperfield of Definition \ref{defnkras}. Recall that $\sign=\{-1,0,1\}$ with hyper-addition given by the ``rule of signs" \eqref{addsign},
and the (classical) multiplication also given by the rule of signs. In this section, we generalize the results proved in \S\ref{hyperequ} for extensions of the hyperfield $\kras$, to hyperfield extensions of $\sign$. In particular, we show that one can recast the hyperaddition in a hyperfield extension of $\sign$ by implementing an {\em order relation} commuting with its conjugates.

\begin{prop}\label{krasnersign}
  Let $R$ be a commutative ring and let $G\subset R^\times$ be a subgroup of its multiplicative group. Assume that $-1\notin G\neq \{1\}$. Then,
  the hyperring $R/G$ contains $\sign$ as a sub-hyperfield if and only if $\{0\}\cup G\cup(-G)$ is an ordered subfield of $R$ with positive part $\{0\}\cup G$.
  \end{prop}
  \proof Let $F=\{0\}\cup G\cup(-G)$. If $(F,G)$ is an ordered field, then $F/G=\sign$ and $R/G$ contains $\sign$ as a sub-hyperfield.  Conversely  one notices that $H=G\cup(-G)$ is a multiplicative subgroup $H\subset R^\times$ and that the hyperring $R/H$ contains $\kras$ as a sub-hyperfield. Thus, by Proposition \ref{krasner1}, $\{0\}\cup G\cup(-G)$ is a subfield $F$ of $R$. This subfield is ordered by the subset $\{0\}\cup G=F_+$. Indeed, from $1+1=1$ in $\sign$ one gets that $G+G=G$ and for $x,y\in F_+$, both $x+y$ and $xy$ are in $F_+$. The statement follows.
\endproof

\begin{prop}\label{addequsign} Let $R$ be a hyperring containing $\sign$ as a sub-hyperring. Then, the following condition
 defines a partial order relation $\leq_R$ on $R$
 \begin{equation}\label{order}
    x \leq_R y~ \Longleftrightarrow~y\in x+1\, \ \text{or}\ \ y=x.
 \end{equation}
\end{prop}
\proof  We show that the relation \eqref{order} is transitive. Assume $x\leq_R y$ and $y\leq_R z$. Then unless one has equality one gets $y\in x+1$ and $z\in y+1$ so that
$z\in (x+1)+1=x+1$ since $1+1=1$. It remains to show that if $x\leq_R y$ and $y\leq_R x$ then $x=y$. If these conditions hold and $x\neq y$ one gets $x\in y+1\subset (x+1)+1=x+1$. Thus $x\in x+1$ but by the reversibility condition $(5)$ on hypergroups one has $1\in x-x$ but $x-x=\{-x,0,x\}$ and one gets that $x=\pm 1$. Similarly $y=\pm 1$, and since $x\neq y$, one of them say $x$ is equal to $1$ and one cannot have $y\in x+1=1$.\endproof

\begin{cor}\label{addequ1sign} Let $R$ be a hyperring containing $\sign$ as a sub-hyperring and let $\leq_R$ be the corresponding partial order relation.  Then
\begin{equation}\label{sxequsign}
    x+1=\{y\geq_R x\,, \ y\neq x\}\qqq x\in R,\ x\neq \pm 1.
\end{equation}
When $R$ is a hyperfield, its additive structure is uniquely determined by the partial order relation $\leq_R$.
\end{cor}

\proof By \eqref{notin}, if $x\neq \pm 1$ one has $x\notin x+1$ and thus using  \eqref{order} one gets \eqref{sxequsign}. This determines the operation $x\mapsto x+1$ for all $x$, including  for $x\in\sign\subset R$.  When $R$ is a hyperfield this determines the addition. \endproof
\begin{cor}\label{noextsign} Any hyperfield extension of $\sign$ is infinite.
\end{cor}

\proof Let $F$ be a hyperfield extension of $\sign$, and $x\in F$, $x\notin \sign\subset F$. Then $(x+1)\cap \sign=\emptyset$, since otherwise using reversibility, one would obtains $x\in \sign$. Let $x_1\in x+1$. Then, one has $x<_F x_1$ and iterating this construction one defines a sequence
$$
x<_F x_1<_F x_2<_F \cdots <_F x_n.
$$
The antisymmetry of the partial order relation shows that the $x_k$ are all distinct.\endproof

\begin{lem}\label{comrelsign} Let $R$ be a hyperring containing $\sign$ as a sub-hyperring and let $\leq_R$ be the corresponding partial order relation.  Then, $\leq_R$ commutes with its conjugates under multiplication by any element $a\in R^\times$.
\end{lem}

\proof Let $T=\leq_R$. One has $T(x)=(x+1)\cup x$ for all $x\in R$. It follows that for the conjugate relation $T^a=aTa^{-1}$ one obtains $T^a(x)=(x+a)\cup x$. Thus
$$
T\circ T^a(x)=(x+a)+1\cup(x+1)\cup(x+a)\cup x=T^a\circ T(x).
$$
\endproof
\begin{prop}\label{lemaddsign} Let  $H$ be an abelian group and let $1\neq \epsilon\in H$ be an element of order two.
Let $S$ be a partial order relation on the set $R=H\cup \{0\}$ such that\vspace{.05in}

$\bullet$~$S(\epsilon)=\{\epsilon,0,1\}$, $S(0)=\{0,1\}$, $S(1)=1$ and
   \begin{equation}\label{rev}
    x\leq_S y~\Longleftrightarrow~ \epsilon y\leq_S \epsilon x
   \end{equation}

$\bullet$~The map $s$ defined by $s(\epsilon)=\{\epsilon,0,1\}$, $s(0)=1$, $s(1)=1$ and
   \begin{equation}\label{ssign}
    s(x)=\{y\in S(x)\,, \ y\neq x\}\qqq x\in R,\ x\notin \{\epsilon,0,1\}
   \end{equation}
  fulfills $s(x)\neq \emptyset$ for all $x$ and  commutes with its conjugates
  for the action of $H$ by multiplication on  $R$.\vspace{.05in}

Then, the hyperoperation
\begin{equation}\label{mainadddefn1}
    x+y:=\begin{cases}  y&\text{if}~x=0\\
xs(yx^{-1})&\text{if}~x\neq 0\end{cases}
\end{equation}
defines a commutative hypergroup law on $R$. With this law as addition, the monoid   $R$ becomes a commutative hyperfield containing $\sign$.
\end{prop}

\proof For $x\in R^\times$, let  $s^x$ be the conjugate of $s$ by multiplication by $x$, \ie
$$
s^x(y):=xs(yx^{-1})\qqq y\in R.
$$
The commutation $s\circ s^x=s^x\circ s$ gives, when applied to $y=0$ and using $s(0)=1$, and $s^x(0)=x$
$$
s(x)=xs(x^{-1}).
$$
Assume that $x\neq 0,y\neq 0$, then
$$
x+y=xs(yx^{-1})=yXs(X^{-1})=ys(X)=y+x\, , \ X=xy^{-1}
$$
The same result holds if $x$ or $y$ is $0$ (if they are both zero one gets $0$, otherwise one gets $0+y=y=y+0$
since $s(0)=1$). Moreover, one has the commutation
$$
s^x\circ s^z=s^z\circ s^x\qqq x,z \in R^\times
$$
which shows that, provided both $x$ and $z$ are non-zero
$$
(x+y)+z=s^z(s^x(y))=s^x(s^z(y))=x+(y+z).
$$
If $x=0$ one has $(x+y)+z=y+z=x+(y+z)$,  and similarly for $z=0$. Thus the addition is associative. The distributivity follows from the homogeneity of \eqref{mainadddefn1}.

 Next, we show that $\forall x  \in R,  \ \exists!~ y(=-x)\,, \ 0\in x+y$. Take $y=\epsilon x$ then, provided $x\neq 0$, one has $x+y=xs(\epsilon)=\{\epsilon x ,0,x\}$ so that $0\in x+y$. We show that $y=\epsilon x$ is the unique solution. For $x\neq 0$ and $0\in x+y$ one has $0\in s(yx^{-1})$. Thus it is enough to show that $0\in s(a)$, $a\neq 0$, implies  $a=\epsilon$. Indeed,  one has $a\leq_S 0$ and thus $0\leq_S \epsilon a$ by \eqref{rev}, thus $\epsilon a=1$,
since $S(0)=\{0,1\}$.

Finally one needs to show that $x\in y+z\Rightarrow z\in x+\epsilon y$. One can assume that $y=1$
using distributivity. We thus need to show that
$$
x\in s(z)\Rightarrow z\in x+\epsilon
$$
In fact, it is enough to show that $\epsilon z\geq_S \epsilon x$ and this holds   by \eqref{rev}. \endproof

\begin{example}{\rm Let $F=U(1)\cup\{0\}$ be the union of the multiplicative group $U(1)$ of complex numbers of modulus one with $\{0\}$. $F$ is, by construction, a multiplicative monoid. For $z,z' \in U(1)$, let $(z,z')\subset U(1)$ be the shortest open interval between $z$ and $z'$. This is well defined if $z'\neq \pm z$. One defines the hyper-addition in $F$ so that $0$ is a neutral element and for $z,z' \in U(1)$ one sets
\begin{equation}\label{uone}
    z+z'=
    \left\{
      \begin{array}{ll}
        z, & \hbox{if} \ z=z'\\
        \{-z,0,z\}, & \hbox{if}\ z=-z'\\
        (z,z'), & \hbox{otherwise.}
      \end{array}
    \right.
\end{equation}
This determines the hyperfield extension of $\sign$: $F=\C/\R_+^\times$. This hyperfield represents the notion of the argument of a complex number. The quotient topology is quasi-compact,  and $0$ is a closed point such that $F$ is its only neighborhood. The subset $U(1)\subset F$ is not closed but the induced topology is the usual topology of $U(1)$.
}\end{example}

\begin{rem}\label{orderskew}{\rm In \S \ref{sectproj} we showed that $\kras$-vector spaces
are projective geometries. Similarly, one can interpret $\sign$-vector spaces
in terms of {\em spherical}  geometries. In the Desarguesian case, any
such geometry is the quotient $V/H^+$ of a left $H$-vector space $V$ over an
ordered skew field $H$ by the {\em positive} part $H^+$ of $H$. It is a double
cover of the projective space $\P(V)=V/H^\times$. More generally, a
$\sign$-vector space $E$ is a double cover of the  $\kras$-vector space
$E\otimes_\sign\kras$. We shall not pursue further this viewpoint in this
paper, but refer to Theorem 28 of Chapter I of \cite{Corsini} as a starting point.
This extended construction is the natural
framework for several results proved in this section. }\end{rem}

\section{Relation with $\B$ and $\F_1$}

By definition, a map $f:R_1\to R_2$ from a hypersemiring $R_1$ to a hypersemiring $R_2$ is a homomorphism when it is a morphism of multiplicative monoids and it fulfills the inclusion
\begin{equation}\label{homo}
    f(x+y)\subset f(x)+f(y)\qqq x,y\in R_1.
\end{equation}
Thus  there is no homomorphism of hypersemirings $f:\Z\to\B$ to the semifield $\B=\{0,1\}$ ($1+1=1$ in $\B$, \cf\cite{Lescot}) such that $f(0)=0, \,f(1)=1$. Indeed $f(-1)$ should be an additive inverse of $1$ in $\B$ which is a contradiction.  On the other hand, the similar map $\sigma:\Z\to\sign$ does extend to a hyperring homomorphism
\begin{equation}\label{hmorp}
\sigma:\Z\to \sign\,, \ \ \sigma(n)={\rm sign}(n)\qqq n\neq 0\,, \ \sigma(0)=0.
\end{equation}
Such map is in fact the unique element of $\Hom(\Z,\sign)$. Moreover, the identity map ${\rm id}:\B\to\sign$ is a hypersemiring homomorphism since $1+1=1$ in $\sign$. Thus one can identify $B$ as the positive part of $\sign$:  $\B=\sign_+$. Notice also that $\kras$ is the quotient of $\sign$ by the subgroup $\{\pm 1\}$. One deduces a canonical epimorphism (absolute value) $\pi:\sign \to \kras$. Thus, by considering the associated geometric spectra (and reversing the arrows), we obtain the following commutative diagram

\begin{gather}
\label{overall}
 \,\hspace{100pt}\raisetag{-47pt} \xymatrix@C=25pt@R=25pt{
 &
    \Spec\sign\ar[dl]\ar[d]&\Spec\kras \ar[l]\ar[dl]\\
 \Spec\Z\ar[d] &
  \Spec\B\ar[dl]& &\\
\Spec\F_1   & &\\
}\hspace{25pt}
\end{gather}

\subsection{The BC-system as $\Z_\sign=\hat\Z\otimes_\Z\sign$}

It follows from what has been explained above that $\Spec\sign$ sits over $\Spec\Z$ and that the map from $\Spec\kras$ to the generic point of $\Spec\Z$ factorizes through $\Spec\sign$. To introduce in this set-up an algebraic geometry  over $\Spec\sign$, it is natural to try to lift $\Spec\Z$ to an object over $\Spec\sign$. This is achieved by considering the spectrum of the tensor product  $\Z_\sign=\hat\Z\otimes_\Z\sign$, using the natural sign homomorphism $\Z\to \sign$ and the embedding of the relative integers in their profinite completion. Notice that, since the non-zero elements of $\sign$ are $\pm 1$, every element of $\hat\Z\otimes_\Z\sign$ belongs to $\hat\Z\otimes_\Z 1$. Thus the hyperring $\Z_\sign$
is, by construction, the quotient of $\hat \Z$ by the equivalence relation
$$
x\sim y~\Longleftrightarrow~\exists n,m\in \N^\times, \ nx=my.
$$
This is precisely the relation that defines the noncommutative space associated to the BC-system. Geometrically, it corresponds to a fibered product given by the commutative diagram

\begin{gather}
\label{BC}
 \,\hspace{60pt}\raisetag{-47pt} \xymatrix@C=25pt@R=25pt{
 &
  \Spec\Z_\sign\ar[dl]\ar[d]& \\
 \Spec\hat\Z\ar[d] &
  \Spec\sign\ar[dl]& \\
\Spec\Z   & \\
}\hspace{25pt}
\end{gather}

Using the morphism $h=\pi\circ\sigma$ of \eqref{hmorp}, one can perform the extension of scalars from $\Z$ to $\kras$. The relation between $-\otimes_\Z\kras$ and $-\otimes_\Z\sign$ is explained by the following result
\begin{prop}\label{extscalar}
Let $R$ be a (commutative) ring containing $\Q$.  Let $R/\Q^\times$ be the hyperring quotient of $R$ by the multiplicative group $\Q^\times$ of $\Q$. Then one has
\begin{equation}\label{extsca}
    R\otimes_\Z\kras=R/\Q^\times.
\end{equation}
Let $R/\Q_+^\times$ be the hyperring quotient of $R$ by the positive multiplicative group $\Q_+^\times$. Then one has
\begin{equation}\label{extscasign}
    R\otimes_\Z\sign=R/\Q_+^\times.
\end{equation}
\end{prop}

\proof   Every element of $R\otimes_\Z\kras$ arises from an element of $R$ in the form $a\otimes 1_{\kras}$. Moreover one has, for $n\in \Z$, $n\neq 0$
$$
n a\otimes 1_{\kras}=a\otimes h(n){1}_{\kras}=a\otimes 1_{\kras}.
$$
This shows that for any non-zero rational number $q\in \Q^\times$ one has
$$
q a\otimes 1_{\kras}= a\otimes 1_{\kras}.
$$
Thus, since $R/\Q^\times$ is a hyperring over $\kras$, by Proposition \ref{krasner2} one gets \eqref{extsca}. The proof of the second statement is similar. \endproof

When $R=\A_{\Q,\,f}$  the ring $\A_{\Q,\,f}=\hat \Z\otimes_\Z\Q$ of finite ad\`eles over $\Q$, Proposition \ref{extscalar} yields the hyperring $\Z_\sign$. Taking $R=\A_\Q$, the ring of ad\`eles over $\Q$, and tensoring by $\kras$ one obtains the hyperring $\H_\Q$ of ad\`ele classes over $\Q$ (\cf Theorem \ref{thmfine} below).
At the level of spectra one obtains
$$
\Spec(\H_\Q) = \Spec(\A_\Q) \times_{\Spec(\Z)} \Spec(\kras)
$$
where $\H_\Q$ is the hyperring of ad\`ele classes over $\Q$. When combined with \eqref{overall}, this construction produces the following (commutative) diagram
\begin{gather}
\label{overall2}
 \,\hspace{40pt}\raisetag{-47pt} \xymatrix@C=25pt@R=25pt{
 &
  \Spec\H_\Q\ar[dl]\ar[d]& \\
 \Spec\A_\Q \ar[d]&
  \Spec\kras\ar[dl]\ar[d]& \\
 \Spec\Z\ar[d] &
  \Spec\B\ar[dl]& \\
\Spec\F_1   & \\
}\hspace{25pt}
\end{gather}

\subsection{The profinite completion $\Z\to\hat\Z$ and ideals}

Let us consider the (compact) topological ring $R=\hat\Z$. Given a closed ideal $J\subset R$, we define
\begin{equation}\label{quasirad}
    \sqrt[\infty]{J}=\{x\in R\mid \lim_{n\to\infty}x^n\subset J\}
\end{equation}
In this definition we are not assuming that $x^n$ converges and we define $\lim_{n\to\infty}x^n$ as the set of limit points of the sequence $x^n$. Thus $x \in \sqrt[\infty]{J}$ means that $x^n\to 0$ in the quotient (compact) topological ring $R/J$.
\begin{lem} \label{basictop} Let $\ell\in \Sigma_\Q$ be a finite place.

$(a)$~For $a=(a_w)\in \hat\Z\sim\prod \Z_p$, the condition $a_\ell=0$ defines a closed ideal $\ffp_\ell\subset \hat\Z$ which is invariant under the equivalence relation induced by the partial action of $\Q^\times$ on $\hat\Z$
by multiplication.

$(b)$~The intersection $\Z\cap \sqrt[\infty]{\ffp_\ell}$ is the prime ideal  $(\ell)\subset\Z$.

\end{lem}

\proof $(a)$~The ideal $\ffp_\ell$ is closed in $\hat\Z\sim\prod \Z_p$  by construction. For any prime $\ell$, the ring $\Z_\ell$ contains $\Z$ and has no zero divisor, thus  $a_\ell=0 \Leftrightarrow n a_\ell=0$ for any non-zero $n\in \N$.

$(b)$~For $a=(a_w)\in \hat\Z\sim\prod \Z_p$, one has $a\in \sqrt[\infty]{\ffp_\ell}$ if and only if the component $a_\ell$ belongs to the maximal ideal $\ell \Z_\ell$. The result follows.\endproof

The relations between the various algebraic structures discussed above are summarized by the following diagram, with $\Z_\kras=\Z_\sign\otimes_\sign \kras=\Z_\sign/\{\pm 1\}$,

\begin{gather}
\label{overall3}
 \,\hspace{50pt}\raisetag{57pt} \xymatrix@C=25pt@R=25pt{
 &
  \Spec\Z_\sign\ar[dl]\ar[d]&\Spec\Z_\kras\ar[d]\ar[l]\ar[r] &\Spec\H_\Q \ar[dl]&\\
 \Spec\hat \Z \ar[d]&
  \Spec\sign\ar[dl]\ar[d]&\Spec\kras \ar[l]\ar[dl]&&\\
 \Spec\Z\ar[d] &
  \Spec\B\ar[dl]& &&\\
\Spec\F_1   & &&\\
}\hspace{25pt}
\end{gather}

\section{Arithmetic of the hyperring $\ads$  of ad\`ele classes} \label{hyper}

The quotient construction of Proposition~\ref{krasner2} applies, in particular, to the ring $R=\A_\K$ of ad\`eles over a global field $\K$, and to the subgroup $\K^\times\subset\A_\K^\times$. One then obtains a new algebraic structure and description of the ad\`ele class space as follows
\begin{thm}\label{thmfine} Let $\K$ be a global field.
The ad\`ele class space $\A_\K/\K^\times$ is a hyperring $\ads$ over $\kras$.  By using the unique morphism $\K\to \kras$ for the extension of scalars one has $\ads=\A_\K\otimes_\K\kras$.
\end{thm}

\proof The fact that $\ads=\A_\K/\K^\times$ is a hyperring follows from the construction of Krasner. This hyperring contains $\kras$ by Proposition~\ref{krasner2}. The identification with $\A_\K\otimes_\K\kras$ follows as in Proposition~\ref{extscalar}.\endproof

This section is devoted to the study of the arithmetic of the hyperring $\H_\K$
of the ad\`ele classes of a global field. In particular we show that,  for
global fields of positive characteristic, the action of the  units
$\H_\K^\times$ on the prime elements of $\H_\K$ corresponds to the action of
the abelianized  Weil group $\cW^{\rm ab}\subset{\rm Gal}(\K^{\rm ab}:\K)$ on
the space ${\rm Val}(\K^{\rm ab})$ of valuations of the maximal abelian
extension $\K^{\rm ab}$ of $\K$ \ie on the space of the (closed) points of the
corresponding projective tower of algebraic curves. More precisely we shall
construct a canonical isomorphism of the groupoid of prime elements of $\H_\K$
with the loop groupoid of the above abelian cover.

\subsection{The space $\Spec(\ads)$ of closed prime ideals of $\ads$}

Given a {\em finite} product of fields $R=\prod_{v\in S}\K_v$, an ideal $J$ in the ring $R$ is of the form
$$
J_Z=\{x=(x_v)\in R\mid x_w=0\qqq w\in Z\}
$$
where $Z\subset S$ is a non-empty subset of $S$. To see this fact one notes that if $x\in J$ and the component $x_v$ does not vanish, then the characteristic function $1_v$ (whose all components are zero except at $v$ where the component is $1$) belongs to $J$ since it is a multiple of $x$. By adding all these $1_v$'s, one gets a generator $p=\sum 1_v$ of $J$.

Let $\K$ be a global field. We endow the ring $\A_\K$ of ad\`eles with its locally compact topology. For any subset $E\subset \Sigma(\K)$ of the set of places of $\K$, one has the convergence
\begin{equation}\label{converge}
\sum_F \, 1_v\to 1_E
\end{equation}
where $F$ runs through the finite subsets of $E$, and $1_E$ is the characteristic function of $E$.

\begin{prop} \label{propidealsads} There is a one to one correspondence between  subsets $Z\subset \Sigma(\K)$ and closed ideals of $\A_\K$ (for the locally compact topology)  given by
\begin{equation}\label{idealsads}
Z\mapsto J_Z=\{x=(x_v)\in \A_\K\mid x_w=0\qqq w\in Z\}.
\end{equation}
\end{prop}

\proof First of all we notice that, by construction, $J_Z$ is a closed ideal of $\A_\K$, for any subset $Z\subset \Sigma(\K)$. Let $J$ be a closed ideal of $\A_\K$. To define the inverse of the map \eqref{idealsads},  let $E\subset \Sigma(\K)$ be the set of places $v$ of $\K$ for which there exists an element of $J$ which does not vanish at $v$. One has $1_v\in J$ for all $v\in E$. Thus, since $J$ is closed one has $1_E\in J$ by \eqref{converge}. The element $1_E$ is a generator of $J$, since for $j\in J$ all components $j_w$ of $j$ vanish for $w\notin E$, so that $j=j1_E$ is a multiple of $1_E$. By taking $Z=E^c$ to be the complement of $E$ in $\Sigma(\K)$, one has $J=J_Z$.
\endproof

Proposition \ref{propidealsads} applies, in particular, in the case $Z = \{w\}$, for  $w\in\Sigma(\K)$ and it gives rise to the closed ideal of the hyperring $\H_\K=\A_\K/\K^\times$
\begin{equation}\label{prideal}
\ffp_w=\{x\in \H_\K\,|\, x_w=0\}.
\end{equation}
Notice that the ideal $\ffp_w$  is well defined since  the condition for an ad\`ele to vanish at a place is invariant under multiplication by elements in $\K^\times$. The set $\ffp_w$ is in fact a {\em prime ideal} in $\H_\K$ whose complement is the multiplicative subset
$$
\ffp_w^c=\{x\in \H_\K\,|\, x_w\neq 0\}.
$$
\begin{prop} \label{propidealsads1} There is a one to one correspondence between the set $\Sigma(\K)$ of places of $\K$ and the set of closed prime ideals of $\H_\K$ (for the quotient topology)  given by
\begin{equation}\label{idealsads1}
\Sigma(\K)\to \Spec(\H_\K),\quad w\mapsto \ffp_w.
\end{equation}
\end{prop}

\proof  The projection $\pi: \A_\K\to \ads$ gives a one to one correspondence for closed prime ideals. Thus, it is enough to prove the statement for the topological ring $\A_\K$. One just needs to show that an ideal of the form $J_Z$ in $\A_\K$ is prime if and only if $Z=\{w\}$ for some place $w\in\Sigma(\K)$. Assume that $Z$ contains two distinct places $w_j$ ($j=1,2$). Then one has $1_{w_j}\notin J_Z$, while the product $1_{w_1}1_{w_2}=0$. Thus $J_Z$ is not a prime ideal of $\A_\K$. Since we have just proved that the $\ffp_w$'s are prime ideals of $\ads$, we thus get the converse.\endproof

\begin{rem}{\rm When viewed as a multiplicative monoid, the ad\`ele class space $\A_\K/\K^\times$ has many more prime ideals than when it is viewed  as a hyperring. This is a consequence of the fact that  {\em any} union of prime ideals in a monoid is still a prime ideal. This statement implies, in particular, that all subsets of the set of places determine a prime ideal.
}\end{rem}

\subsection{Functions on $\Spec(\H_\Q)$}\label{functions}

In algebraic geometry one defines a function on a scheme $X$, viewed as a covariant functor $\underline X:\An\to\Se$, as  a morphism of functors $f:\underline X\to \cD$ to the (functor) affine line $\cD=\mathfrak{spec}(\Z[T])$ (whose geometric scheme is $\Spec(\Z[T])$, \cf \cite{announc3}). For $X=\Spec(R)$, where $R$ is an object of $\An$ (\ie a commutative ring with unit), one derives a natural identification of functions on $X$ with elements of the ring $R$
\begin{equation}\label{ztidd}
   \Hom_\An(\Z[T],R)\simeq R.
\end{equation}
 In the category of hyperrings, the identification \eqref{ztidd} no longer holds in general as easily follows from Proposition \ref{ex1}. Indeed, $\kras$ has only two elements while $\Hom_\han(\Z[T],\kras)\simeq \Spec(\Z[T])$ is countably infinite. The next Theorem describes the functions on $\Spec(\ads)$, for $\K=\Q$.
 \begin{thm} \label{functads} Let $\H_\Q$ be the hyperring of ad\`ele classes over $\Q$, and let $\rho\in \Hom_\han(\Z[T],\H_\Q).$ Then, either $\rho=\xi_a$
 \begin{equation}\label{lift}
    \xi_a(P(T))=P(a)\Q^\times \in \H_\Q\qqq P\in \Z[T]
 \end{equation}
 for a unique ad\`ele $a\in \A_\Q$, or $\rho$ factorizes through $\Q[e_Z]/\Q^\times$, where $e_Z$  is the idempotent of $\A_\Q$
 associated to a subset $Z\subset \Sigma_\Q$.
 \end{thm}
 \proof  Assume first that the range $\rho(\Z[T])$ is contained in $L\cup \{0\}$ where $L$ is a line of the projective space  $\H_\Q\backslash \{0\}$. Let $\pi:\A_\Q\to \H_\Q$ be the projection. The two dimensional subspace $E=\pi^{-1}(L\cup \{0\})$ of the $\Q$-vector space $\A_\Q$ contains $1$ since $\rho(1)=1$. Unless  $\rho(\Z[T])=\kras$,
 the line $L$ is generated by $1$ and $\xi\in\rho(\Z[T])$, $\xi\neq 1$.
 Let $x\in E$ with $\pi(x)=\xi$. Since $\xi^2\in\rho(\Z[T])$ one has $x^2\in E$
 and $x^2=ax + b$ for some $a,b\in\Q$. As in remark \ref{sub2}, one can assume that $x^2=N$ for a square free integer $N$.
 But the equation $y^2=N$ has no solution in $\A_\Q$ except for $N=1$. It follows that $\Q[x]\subset \A_\Q$ is
  a two dimensional subalgebra over $\Q$  of the form
 $\Q[e_Z]$ where $e_Z$ is the idempotent of $\A_\Q$
 associated to a subset $Z\subset \Sigma_\Q$. Thus $\rho$ factorizes through $\Q[e_Z]/\Q^\times$.
 We can thus assume now that   $\rho(\Z[T])\backslash \{0\}$ is not contained in a line $L$ of the projective space  $\H_\Q\backslash \{0\}$. The restriction of $\rho$ to $\Z\subset \Z[T]$ is a morphism from $\Z$ to $\kras$ and its kernel is a prime ideal $\ffp\subset \Z$. If $\ffp\neq \{0\}$, one has $\ffp=p\Z$ for a prime number $p$. Then
 $$
 \rho\left(\sum (pa_n)T^n\right)=\rho(p)\rho\left(\sum a_nT^n\right)=0\qqq a_k\in \Z
 $$
 and the inclusion $\rho(x+y)\subset\rho(x)+\rho(y)$ shows that $\rho(P(T))$ only depends upon the class
 of $P(T)$ in $\F_p[T]$.  Since $\rho(\F_p^\times)=1$ one gets a morphism\footnote{Note that this holds even for $p=2$
 even though $\F_2[T]/\F_2^\times$ is not a $\kras$-vector space.},
 in the sense of \cite{Faure}, from the projective space $(\F_p[T]/\F_p^\times)\backslash \{0\}$ to the projective space  $\H_\Q\backslash \{0\}$. Since $\rho(\Z[T])\backslash \{0\}$ is not contained in a line $L$ of the projective space  $\H_\Q\backslash \{0\}$, this morphism is non-degenerate. By \opcit Theorem 5.4.1, (\cf also \cite{Faure1} Theorem 3.1) there exists a semi-linear map inducing this morphism but this gives a contradiction since there is no field homomorphism from $\F_p$ to $\Q$. Thus one has $\ffp=\{0\}$ and $\rho(n)=1$ for all $n\in \Z\backslash \{0\}$. One can then extend $\rho$ to a morphism
 $$
  \rho'\,:\, \Q[T]\to \H_\Q\,, \ \  \rho'(P(T))=\rho(nP(T))\qqq n \neq 0, \ nP(T)\in \Z[T].
 $$
By Corollary \ref{main1} one then gets a unique ring homomorphism $\tilde\rho: \Q[T]\to \A_\Q$ which lifts $\rho'$. This gives a unique ad\`ele $a\in \A_\Q$ such that \eqref{lift} holds.\endproof

The above result shows that there are two different types of ``functions" on $\Spec(\H_\Q)$: functions corresponding to ad\`eles (which recover the algebraic information of the ring $\A_\Q$) and functions  factorizing through $\Q[e_Z]/\Q^\times$. These latter functions should be thought of as {\em ``two-valued"} functions, in analogy with the case of continuous functions  on a compact space $X$. Indeed, the range of $f\in C(X,\R)$ has two elements if and only if the subalgebra of $C(X,\R)$ generated by $f$ is of the form $\R[e]$ for some idempotent $e\in C(X,\R)$. In the above case of  $\Spec(\H_\Q)$ the subset $Z\subset \Sigma_\Q\simeq \Spec(\H_\Q)$ and its complement specify the partition of $\Spec(\H_\Q)$ corresponding to the two values of $\rho$. Once this partition is given, the remaining freedom is in the set  $\Hom_\han(\Z[T],(\Q\oplus\Q)/\Q^\times)$.
We shall not attempt to describe explicitly this set here, but refer to Remark \ref{sub2} to show that it contains many elements.\vspace{.05in}

 Let $\cH$ be a commutative ring, and let $\Delta:\cH\to \cH\otimes_\Z\cH$ be a coproduct. Given two ring homomorphisms $\rho_j: \cH\to R$ ($j=1,2$) to a commutative ring $R$, the composition
$\rho=(\rho_1\otimes \rho_2)\circ \Delta$ defines a homomorphism $\rho: \cH\to R$. When $R$ is a hyperring,
one introduces the following notion
\begin{defn} Let $(\cH,\Delta)$ be a commutative ring with a coproduct and let $R$ be a hyperring. Let $\rho_j\in \Hom_\han(\cH,R)$, $j=1,2$. One defines  $\rho_1\star_\Delta \rho_2$ to be the set of $\rho\in \Hom_\han(\cH,R)$ such that for any $x\in \cH$ and any decomposition $\Delta(x)=\sum x_{(1)}\otimes x_{(2)}$ one has
\begin{equation}\label{equhopf}
    \rho(x)\in \sum \rho_1(x_{(1)})\rho_2(x_{(2)})  \,.
\end{equation}
\end{defn}
 In genera,l $\rho_1\star_\Delta \rho_2$ can be empty or it may contain several elements. When $\rho_1\star_\Delta \rho_2=\{\rho\}$ is made by a single element we simply write  $\rho_1\star_\Delta \rho_2= \rho $.\vspace{.05in}

When $\cH=\Z[T]$, $\Delta^+(T)=T\otimes 1+1\otimes T$ and $\Delta^\times(T)=T\otimes T$,  this construction
allows  one to recover  the algebraic structure of the ring of ad\`eles, in terms of functions on $\Spec(\H_\Q)$ (\cf Theorem \ref{functads}).
\begin{prop}\label{Hopf} Let $\rho_j=\xi_{a_j}\in \Hom_\han(\Z[T],\H_\Q)$ ($j=1,2$), be the homomorphisms uniquely associated to $a_j\in \A_\Q$ by \eqref{lift}. Assume that monomials of degree $\leq 2$ in $a_j$ are linearly independent over $\Q$. Then one has
 \begin{equation}\label{coprodfine}
 \rho_1\star_{\Delta^+} \rho_2=\xi_{a_1+a_2}\,, \ \ \rho_1\star_{\Delta^\times} \rho_2=\xi_{a_1a_2}.
 \end{equation}
\end{prop}

\proof Let $\tilde \rho_j\in \Hom(\Z[T],\A_\Q)$, $\tilde \rho_j(P(T))=P(a_j)$ be the lift of $\rho_j$. Then  $\rho^+=(\tilde\rho_1\otimes \tilde\rho_2)\circ \Delta^+$ fulfills the equation
$$
\rho^+(x)=\sum \tilde\rho_1(x_{(1)})\tilde\rho_2(x_{(2)})
$$
for any decomposition $\Delta^+(x)=\sum x_{(1)}\otimes x_{(2)}$. Thus, since $\xi_{a_1+a_2}=\pi\circ  \rho^+$ (with $\pi:\A_\Q\to \H_\Q$ the projection) one gets, using \eqref{longsum}, that $\xi_{a_1+a_2}\in \rho_1\star_{\Delta^+} \rho_2$. In a similar manner one obtains $\xi_{a_1a_2}\in \rho_1\star_{\Delta^\times} \rho_2$. It remains to show that they are the only solutions. We do it first for $\Delta^\times$. Let $\rho\in \rho_1\star_{\Delta^\times} \rho_2$. Since $\Delta^\times(T)=T\otimes T$, \eqref{equhopf} gives
$\rho(T)=a_1a_2\in \H_\Q$. Similarly $\rho(T^2)=a_1^2a_2^2\in \H_\Q$. Thus since $1,a_1a_2, a_1^2a_2^2$ are linearly independent over $\Q$, the range of $\rho$ is of $\kras$-dimension $\geq 3$ and by Theorem \ref{functads} there exists $a\in \A_\Q$ such that $\rho=\xi_a$. Moreover $a=\lambda a_1a_2$ for some $\lambda\in\Q^\times$ and it remains to show that $\lambda=1$. One has
$$
\Delta^\times(1+T)=(1-T)\otimes (1- T)+T\otimes 1+1\otimes T
$$
Thus \eqref{equhopf} shows that $1+a$ belongs to
$
(1-a_1)(1-a_2)\Q^\times+ a_1 \Q^\times+a_2\Q^\times.
$
But by $\Q$-linear independence the only element of this set which is of the form $1+\lambda a_1 a_2$ is $1+a_1a_2$ which implies that $\lambda=1$ and thus that $\rho=\xi_{a_1a_2}$. Let now $\rho\in \rho_1\star_{\Delta^+} \rho_2$. One has
$$
\Delta^+(T)=T\otimes 1+1\otimes T=(1+T)\otimes (1+T)-1\otimes 1-T\otimes T
$$
which implies by \eqref{equhopf} that $\rho(T)=\lambda(a_1+a_2)$ for some $\lambda\in \Q^\times$. Since $\rho(T^2)\in a_1^2\Q^\times+ a_1a_2 \Q^\times+a_2^2\Q^\times$ the range of $\rho$ is of $\kras$-dimension $\geq 3$ and by Theorem \ref{functads} there exists $a\in \A_\Q$ such that $\rho=\xi_a$. One has $a=\lambda(a_1+a_2)$ and to show that $\lambda=1$ one proceeds as above using
$$
\Delta^+(1+T)=(1+T)\otimes (1+T)-T\otimes T\,.
$$
\endproof

\subsection{The groupoid $P(\ads)$ of prime elements of $\ads$}\label{space}

The notion of {\em principal} prime ideal in a hyperring is related to the following notion of prime element
\begin{defn}\label{prime} In a  hyperring $R$,  an element $a\in R$ is said to be {\em prime} if the ideal $aR$ is a prime ideal.
\end{defn}
We let $P(\ads)$ be the set of prime elements of the hyperring $\ads=\A_\K/\K^\times$, for $\K$ a global field.
\begin{thm}\label{classprime}~$1)$~Any  principal prime ideal of $\ads$ is equal to $\ffp_w = a\H_\K$ for a  place $w\in\Sigma(\K)$ uniquely determined by $a\in\H_\K$.\vspace{.05in}

$2)$~The group $C_\K=\A_\K^\times/\K^\times$ acts transitively on the generators of the principal prime ideal $\ffp_w$.\vspace{.05in}

$3)$~The isotropy subgroup of any generator of the  prime ideal $\ffp_w$ is $\K_w^\times\subset C_\K$.
\end{thm}
\proof $1)$~ Let $\ffp=a\ads$ be a prime principal ideal in $\ads$. We consider the support of $a$ \ie the set $
S=\{v\in \Sigma(\K)\,|\, a_v\neq 0\}$.
  We shall prove that the characteristic function $1_S$ generates the same ideal as $a$, \ie $a\ads=1_S \ads$, where $1_S\in \ads$ is the class of the ad\`ele $\alpha= (\alpha_v)$, with $\alpha_v=1$ for $v\in S$ and $\alpha_v=0$ otherwise.
For each $v\in \Sigma(\K)$, let $\cO^\times_v$ be the multiplicative group $
\cO^\times_v=\{x\in \K_v\,:\,|x|_v=1\}$
of elements in $\K_v$ of norm $1$. When the place $v$ is non-archimedean this is the group of invertible elements of the local ring $\cO_v$.
We let $a=(a_v)$ be an ad\`ele in the given class and first show that there exists a finite subset $F\subset \Sigma(\K)$ such that
\begin{equation}\label{contra}
 a_v\in \cO^\times_v\qqq v\in S \,, \ v\notin F.
\end{equation}
Otherwise, there would exist an infinite subset $Y\subset S $ such that
 $$
 |a_v|_v<1\qqq v\in Y.
 $$
 Let then $Y'$ be an infinite   subset of $Y$ whose complement in $Y$ is infinite. Consider the ad\`eles $y,z\in \A_\K$ defined by
 $$
 y_v=\left\{
      \begin{array}{ll}
        1, &   v\notin Y' \\
        a_v, & v\in Y'
      \end{array}
    \right.
    \,, \ \
    z_v=\left\{
      \begin{array}{ll}
        a_v, &   v\notin Y' \\
        1, & v\in Y'.
      \end{array}
    \right.
 $$
 By construction $yz=a$. The same equality holds in $\ads$. Since the ideal $\ffp=a\ads$ is prime, its complement in $\ads$ is multiplicative and thus  $y\in \ffp$ or (and) $z\in \ffp$. However $y\notin a\A_\K$ since $|y_v|_v> |a_v|_v$ on the complement of $Y'$ in $Y$ which is an infinite set of places. Similarly $z\notin a\A_\K$ since $|z_v|_v> |a_v|_v$ on $Y'$. Thus one gets a contradiction and this proves \eqref{contra}. In fact one may assume, without changing the principal ideal $\ffp=a\ads$, that
 \begin{equation}\label{invert}
 a_v\in \cO^\times_v\qqq v\in S
\end{equation}
Since the ideal $\ffp=a\ads$ is non-trivial the complement $Z$ of $S $ in $\Sigma(\K)$ is non-empty.
Assume that $Z$ contains two places $v_1\neq v_2$. Let $1_{v}$ be the (class of the) ad\`ele all of whose components vanish except at the place $v$ where its component is $1\in\K_v$. Then one has $1_{v_j}\notin \ffp=a\ads$, but the product
 $1_{v_1}1_{v_2}=0\in \ffp=a\ads$ which  contradicts the fact that $\ffp=a\ads$ is prime. This shows that $Z=\{w\}$ for some $w\in\Sigma(\K)$ and that $a\ads=\ffp_w$ using \eqref{invert}.

$2)$~ Let $a,b\in \ads$ be two generators of the ideal $\ffp_w$. Let $\alpha,\beta\in \A_\K$ be ad\`eles in the classes of $a$ and $b$ respectively. Then by  \eqref{contra}, the equality
$$
j_v=\beta_v/\alpha_v\qqq v\neq w\,, \ \ j_w=1
$$
defines an id\`ele such that $j\alpha=\beta$. This shows that the group $C_\K$ acts transitively on the generators of $\ffp_w$.

$3)$~ Let $a\in \ads$ be a generator of the principal ideal $\ffp_w$ and let $\alpha$ be an ad\`ele in the class of $a$. For $g\in C_\K$ the equality $ga=a$ in $\ads$ means that for $j$ an id\`ele in the class of $g$, there exists $q\in \K^\times$ such that $j\alpha=q\alpha$. In other words one has $q^{-1}j\alpha=\alpha$. Since all components $\alpha_v$ are non-zero except at $v=w$ one thus gets that all components of $q^{-1}j$ are equal to $1$ except at $w$. The component $j_w$ can be arbitrary and thus, the isotropy subgroup of any generator of $\ffp_w$ is $\K_w^\times\subset C_\K$.
\endproof

On $P(\ads)$ we define a groupoid law given by multiplication. More precisely,

\begin{prop}\label{groupoid} Let $\K$ be a global field and $s:P(\ads)\to \Sigma_\K$ the map which associates to a prime element of $\ads$ the principal prime ideal of $\ads$ it generates. Then $P(\ads)$ with range and source maps equal to $s$ and partial product given by the product in the hyperring $\ads$,
is a groupoid.
\end{prop}

\proof Since the source and range maps coincide, one needs simply to show that each fiber $s^{-1}(v)$ is a group. For each place $v\in \Sigma_\K$, there is a unique generator $p_v$ of the prime ideal $\ffp_v$ which fulfills $p_v^2=p_v$. It is given by the class of the characteristic function $1_S$ where $S$ is the complement of $v$ in $\Sigma_\K$. Any other element of $s^{-1}(v)$ is, by Theorem \ref{classprime}, of the form $\gamma=up_v$ where
$u\in C_\K/\K^\times$ is uniquely determined. The product in $s^{-1}(v)$ corresponds to the product in the group $C_\K/\K^\times$.\endproof

Note that the product $p_1p_2$ of two prime elements is  a prime element only when these factors generate the same ideal.

\subsection{The groupoids $\Pi_1^{\rm ab}(X)'$ and $P(\ads)$ in characteristic $p\neq 0$}

Let $\K$ be  a global field of characteristic $p> 0$ \ie a function field over a constant field $\F_q\subset \K$. We fix a separable closure $\bar \K$ of $\K$ and let
 $\K^{\rm ab}\subset \bar \K$ be the maximal abelian extension of $\K$. Let $\bar \F_q\subset \bar \K$ be the algebraic closure of  $\F_q$. We denote by $\cW^{\rm ab}\subset {\rm Gal}( \K^{\rm ab}:\K)$  the abelianized  Weil group, \ie the subgroup of elements of ${\rm Gal}( \K^{\rm ab}:\K)$ whose restriction to $\bar \F_q$ is an integral power of the Frobenius.

Let ${\rm Val}(\K^{\rm ab})$ be the space of all  valuations of $\K^{\rm ab}$. By restriction to $\K\subset \K^{\rm ab}$ one obtains a natural map
\begin{equation}\label{mapp}
    p\;:\; {\rm Val}(\K^{\rm ab})\to \Sigma_\K\,, \ \ p(v)=v|_\K.
\end{equation}
By construction, the action of ${\rm Gal}( \K^{\rm ab}:\K)$ on ${\rm Val}(\K^{\rm ab})$ preserves the map $p$.

\begin{prop}\label{val}Let $w\in \Sigma_\K$.\vspace{.05in}

    $(1)$~The abelianized  Weil group $\cW^{\rm ab}$ acts transitively on the fiber $p^{-1}(w)$ of $p$.\vspace{.05in}

  $(2)$~The isotropy subgroup of an element in the fiber $p^{-1}(w)$ coincides with the abelianized  local Weil
  group $\cW^{\rm ab}_w\subset \cW^{\rm ab}$.
\end{prop}

\proof This follows from standard results of class field theory but we give the detailed proof for completeness. Let $v\in {\rm Val}(\K^{\rm ab})$ with $p(v)=w$. By construction, the completion $\K^{\rm ab}_v$ of $\K^{\rm ab}$ at $v$  contains the local field $\K_w$ completion of $\K$ at $w$. The subfield $\K_w \vee \K^{\rm ab}$ of $\K^{\rm ab}_v$ generated by $\K^{\rm ab}$ and $\K_w$ coincides with the maximal abelian extension $\K_w^{\rm ab}$ of $\K_w$. One has the translation isomorphism (\cf \cite{Bourbaki} Theorem V, A.V.71)
\begin{equation}\label{galiso}
   {\rm Gal}((\K_w \vee \K^{\rm ab}):\K_w) \cong {\rm Gal}( \K^{\rm ab}:(\K_w\cap \K^{\rm ab}))\subset {\rm Gal}( \K^{\rm ab}:\K)
\end{equation}
obtained by restricting an automorphism to $\K^{\rm ab}$.

The subgroup  ${\rm Gal}( \K^{\rm ab}:(\K_w\cap  \K^{\rm ab}))\subset {\rm Gal}( \K^{\rm ab}:\K)$ is the isotropy subgroup $\Gamma_v$ of the valuation $v$ \ie  an element $g\in {\rm Gal}( \K^{\rm ab}:\K)$ fixes $v$ if and only if $g$ fixes pointwise the subfield $\K_w\cap  \K^{\rm ab}\subset  \K^{\rm ab}$. Indeed, if $g$ fixes $v$ it extends uniquely by continuity to an automorphism of $\K^{\rm ab}_v$. This automorphism is the identity on $\K$ and hence also  on the completion $\K_w$ of $\K$ at $w$ and thus on  $\K_w\cap \K^{\rm ab}$. Next, let $g\in {\rm Gal}( \K^{\rm ab}:\K)$ be the identity on $\K_w\cap  \K^{\rm ab}$. The fact that $g$ fixes $v$  follows from \eqref{galiso}.
  Indeed, this shows that any element
 $g\in {\rm Gal}( \K^{\rm ab}:(\K_w\cap \K^{\rm ab}))$ is the restriction of an automorphism in ${\rm Gal}((\K_w \vee \K^{\rm ab}):\K_w)$ and preserves $v$ since the valuation $w$ of the local field $\K_w$ extends uniquely to finite algebraic extensions of $\K_w$, and thus to $\K_w \vee \K^{\rm ab}$, by uniqueness of the maximal compact subring.

 $(1)$~Let us check that the abelianized  Weil group $\cW^{\rm ab}$ acts transitively on the valuations in the set $p^{-1}(w)$. The Galois group ${\rm Gal}( \K^{\rm ab}:\K)$  acts transitively on $p^{-1}(w)$. Indeed
 the space of valuations extending $w$ is by construction the
projective limit of the finite sets of valuations extending
$w$ over finite algebraic extensions of $\K$. The Galois group
${\rm Gal}( \K^{\rm ab}:\K)$ is a compact profinite group which acts transitively on the finite sets of valuations extending
$w$ over finite algebraic Galois extensions of $\K$
 (\cite{Rosen}, \S~9 Proposition 9.2). Thus it acts
transitively on the fiber $p^{-1}(w)$.   It remains to show that the transitivity of the action continues to hold for $\cW^{\rm ab}\subset {\rm Gal}( \K^{\rm ab}:\K)$.  It is enough to show that the orbit $\cW^{\rm ab} v$ of a place  $v\in p^{-1}(w)$ is the same as its orbit under ${\rm Gal}( \K^{\rm ab}:\K)$. This is a consequence of the co-compactness of the isotropy subgroup $\Gamma_v\cap \cW^{\rm ab}\subset \cW^{\rm ab}$ but it is worthwhile to describe what happens in more details. In the completion process from $\K$ to  $\K_w$, the maximal finite subfield (constant field) passes from $\F_q$  to a finite extension $\F_{q^\ell}$.  Let $\K_w^{\rm un}\subset \K_w^{\rm ab}$ be the largest unramified extension of $\K_w$ inside $\K^{\rm ab}$. It is obtained by adjoining to $\K_w$ all roots of unity of order prime to $p$ which are not already in the constant field $F_{q^\ell}$ of $\K_w$.  One has (\cite{Weil} VI, \cite{Tate} Chapter VII)
$$
 {\rm Gal}(\K_w^{\rm ab}:\K_w^{\rm un}) \cong {\rm Gal}( \K^{\rm ab}:(\K_w^{\rm un}\cap \K^{\rm ab})) \subset {\rm Gal}( \K^{\rm ab}:(\K_w\cap \K^{\rm ab}))
$$
 The extension
$\K_w^{\rm un}\cap  \K^{\rm ab}$ contains $\bar\F_q\otimes_{\F_q}\K$. This determines the following diagram of inclusions of fields
\begin{gather}
\label{functCFTmap}
 \,\hspace{25pt}\raisetag{-47pt} \xymatrix@C=25pt@R=25pt{
\K_w^{\rm ab}= \K_w \vee \K^{\rm ab}\subset (\K^{\rm ab})_v &
 \K^{\rm ab} \ar[l]^-{ }& \\
\K_w^{\rm un} \ar[u]^-{ }  & \K_w^{\rm un}\cap \K^{\rm ab}\ar[u]^{ }\ar[l]^{  }&\bar\F_q\otimes_{\F_q}\K\ar[l]^{  }\\
\K_w \ar[u]^-{ }  & \K_w\cap \K^{\rm ab}\ar[u]^{ }\ar[l]_{  }&\K\ar[u]^{ }\ar[l]^{  }\nonumber \\
}\hspace{140pt}
\end{gather}

The topological generator of ${\rm Gal}(\K_w^{\rm un}:\K_w)$ induces the $\ell$-th power $\sigma^\ell$ of the Frobenius automorphism on $\K'$ and the same holds for the topological generator of ${\rm Gal}((\K_w^{\rm un}\cap \K^{\rm ab}):(\K_w\cap \K^{\rm ab}))$. The abelianized  Weil group $\cW^{\rm ab}\subset {\rm Gal}(\K^{\rm ab}:\K)$ is defined by
$$
\cW^{\rm ab}=\rho^{-1}\{\sigma^\Z\}\,, \ \rho\,:\, {\rm Gal}(\K^{\rm ab}:\K)\to {\rm Gal}(\K':\K)\,, \  \sigma^\Z\subset {\rm Gal}(\K':\K)
$$
Thus, the statement that the group $\cW^{\rm ab}\subset {\rm Gal}(\K^{\rm ab}:\K)$ acts transitively on the fiber $p^{-1}(w)$ is equivalent to the fact that the dense  subgroup
 $\Z\subset\hat\Z$ acts transitively in the finite space $\hat\Z/\ell\hat\Z$.

 $(2)$~follows from $(1)$ and the remarks made at the beginning of the proof.\endproof

We now implement the geometric language.  Given an extension $E$ of $\bar \F_q$ of transcendence degree $1$, it is a well-known fact that the  space of valuations of $E$,  ${\rm Val}(E)$, coincides with the set of (closed) points of the unique projective nonsingular algebraic curve with function field $E$. Moreover, one also knows  (\cf \cite{Hart} Corollary 6.12) that the category of nonsingular projective algebraic curves and dominant morphisms is equivalent to the category of function fields of dimension one over $\bar \F_q$. Given a global field $\K$ of positive characteristic $p>1$ with constant field $\F_q$, one knows that the maximal abelian extension $\K^{\rm ab}$ of $\K$ is an inductive limit of extensions $E$  of $\bar \F_q$ of transcendence degree $1$.
Thus the space ${\rm Val}(\K^{\rm ab})$ of valuations of $\K^{\rm ab}$, endowed with the action of the abelianized  Weil group $\cW^{\rm ab}\subset{\rm Gal}(\K^{\rm ab}:\K)$, inherits the structure of a projective limit of projective nonsingular curves. This construction determines the maximal abelian cover $\pi:X^{\rm ab}\to X$ of the non singular projective curve $X$ over $\F_q$ with function field  $\K$.\vspace{.05in}

Let $\pi:\tilde X\to X$ be a Galois covering of $X$ with Galois group $W$. The fundamental groupoid of $\pi$ is by definition the quotient $\Pi_1=(\tilde X\times \tilde X)/W$ of $\tilde X\times \tilde X$ by the diagonal action of $W$ on the self-product. The (canonical) range and source maps: $r$ and $s$ are defined by the two projections
\begin{equation}\label{rs}
    r(\tilde x,\tilde y)=x\,, \ s(\tilde x,\tilde y)=y.
\end{equation}
Let us consider the subgroupoid of {\em loops} \ie
\begin{equation}\label{sub}
    \Pi_1'=\{\gamma\in \Pi_1\mid r(\gamma)=s(\gamma)\}.
\end{equation}
 Each fiber of the natural projection $r=s:\Pi_1'\to X$ is a group. Moreover, if $W$ is an abelian group  one defines the following natural action of $W$ on $\Pi_1'$
\begin{equation}\label{actionW}
    w\cdot (\tilde x,\tilde y)=(w\tilde x,\tilde y)=(\tilde x,w^{-1}\tilde y).
\end{equation}
We apply these results to the maximal abelian cover $\pi:X^{\rm ab}\to X$ of the non singular projective curve $X$ over $\F_q$ with function field $\K$.

We view $X$ as a scheme over $\F_q$. In this case, we let $W=\cW^{\rm ab}\subset{\rm Gal}(\K^{\rm ab}:\K)$  be the abelianized Weil group. We let $\Pi_1^{\rm ab}(X)$ be the fundamental groupoid of this maximal abelian cover and $\Pi_1^{\rm ab}(X)'\subset \Pi_1^{\rm ab}(X)$ the loop groupoid. Since the two projections from $X^{\rm ab}\times X^{\rm ab}$ to $X$ are $W$-invariant, $\Pi_1^{\rm ab}(X)'$ is the quotient of the fibered product $X^{\rm ab}\times_X X^{\rm ab}$ by the diagonal action of $W$. We identify the closed points of $X^{\rm ab}\times_X X^{\rm ab}$ with pairs of valuations of $\K^{\rm ab}$ whose restrictions to $\K$ are the same.\vspace{.05in}

We obtain the following refinement of Proposition 8.13 of
\cite{CCM2}.

\begin{thm}\label{ccm2prop} Let $\K$ be a  global field of characteristic $p\neq 0$, and let $X$ be the corresponding non-singular projective algebraic curve over $\F_q$.\vspace{.05in}

 $\bullet$~The loop groupoid $\Pi_1^{\rm ab}(X)'\subset \Pi_1^{\rm ab}(X)$ is canonically isomorphic to the groupoid $P(\ads)$ of prime elements of the hyperring $\ads=\A_\K/\K^\times$.\vspace{.05in}

    $\bullet$~The above isomorphism $\Pi_1^{\rm ab}(X)'\simeq P(\ads)$  is equivariant for the action of $W$ on $\Pi_1^{\rm ab}(X)'$  and the action of the units $\ads^\times=C_\K$ on prime elements by multiplication.
 \end{thm}

\proof Under the class-field theory isomorphism $W=\cW^{\rm ab}\sim C_\K$, the local Weil group at a place $w\in\Sigma_\K$ corresponds to the subgroup $\K_w^\times\subset C_\K$. By applying Proposition \ref{val}, this shows that given two elements  $v_j\in{\rm Val}(\K^{\rm ab})$ above the same place $w\in \Sigma_\K$,  there exists a unique element $\gamma(v_1,v_2)\in C_\K/\K_w^\times$ such that (under the class field theory isomorphism)
\begin{equation}\label{move}
    \gamma(v_1,v_2)(v_2)=v_1.
\end{equation}
For  a place $v\in \Sigma_\K$ we let $p_v\in P(\ads)$ be the unique idempotent element (\ie $p^2_v=p_v$)
which generates the ideal $\ffp_v$. We define the map (\cf\eqref{mapp})
\begin{equation}\label{groupoidmap}
    \varphi: \Pi_1^{\rm ab}(X)'\to P(\ads)\,, \  \varphi(v_1,v_2)=\gamma(v_1,v_2)p_w \qqq v_j\in p^{-1}(w).
\end{equation}
The map $\varphi$ is well defined since by  Theorem \ref{classprime}  the  isotropy subgroup of points above $w$ in $P(\ads)$ is $\K_w^\times$ and  one has $\gamma(uv_1,uv_2)=\gamma(v_1,v_2)$ for all $u\in \cW^{\rm ab}\sim C_\K$. One also checks the equivariance
\begin{equation}\label{equiv}
    \varphi(u\cdot \alpha)=u\varphi(\alpha) \qqq u \in \cW^{\rm ab}\sim C_\K.
\end{equation}
Finally, the equality
\begin{equation}\label{complaw}
   \gamma(v_1,v_2)\gamma(v_2,v_3)=\gamma(v_1,v_3)
\end{equation}
together with $ap_vbp_v=abp_v$ show that the map $\varphi$ is a morphism of groupoids which is also bijective
over each place in $\Sigma_\K$, by Proposition \ref{val} and  Theorem \ref{classprime}. Thus $\varphi$ is an isomorphism.
\endproof

\end{document}